\documentclass[a4paper, 12pt, draft, reqno]{amsart}
\usepackage[final]{graphicx}
\usepackage{amsmath}
\usepackage{amssymb}
\usepackage{latexsym}
\usepackage[cmtip,matrix,arrow]{xy}
\UseComputerModernTips
\CompileMatrices

\DeclareMathOperator{\sgn}{sgn}
\DeclareMathOperator{\swap}{s}
\newcommand{\ZZ}{{\mathbb Z}}
\newcommand{\QQ}{{\mathbb Q}}
\newcommand{\RR}{{\mathbb R}}
\newcommand{\NN}{{\mathbb N}}
\newcommand{\CC}{{\mathbb C}}
\newcommand{\FF}{{\mathbb F}}

\DeclareMathOperator{\Hom}{Hom}

\DeclareMathOperator{\im}{im}
\DeclareMathOperator{\hd}{hd}

\DeclareMathOperator{\res}{res}
\DeclareMathOperator{\add}{add}
\DeclareMathOperator{\remo}{rem}
\DeclareMathOperator{\supp}{supp}
\DeclareMathOperator{\abv}{abv}
\DeclareMathOperator{\lft}{lft}

\newcommand{\otherwise}{\mbox{\rm otherwise}}
\newcommand{\wif}{\mbox{\rm if }}
\newcommand{\wand}{\mbox{\rm and }}

\DeclareMathOperator{\ind}{ind}

\newcommand{\too}{\longrightarrow}

\begin{document}
\theoremstyle{plain}
\newtheorem{thm}{Theorem}[section]
\newtheorem{prop}[thm]{Proposition}
\newtheorem{lem}[thm]{Lemma}
\newtheorem{cor}[thm]{Corollary}
\newtheorem{conj}[thm]{Conjecture}
\newtheorem{claim}[thm]{Claim}
\theoremstyle{definition}
\newtheorem{rem}[thm]{Remark}
\newtheorem{ass}[thm]{Assumption}
\newtheorem{defn}[thm]{Definition}
\newtheorem{example}[thm]{Example}

\setlength{\parskip}{1ex}

\title[The blocks of the walled Brauer algebra]{On the blocks of the
walled  Brauer algebra}
 \author{Anton Cox}
 \email{A.G.Cox@city.ac.uk, M.Devisscher@city.ac.uk}
 \author{Maud De Visscher} 
\address{Centre for Mathematical Science\\
 City University\\
 Northampton Square\\ 
 London\\ 
 EC1V 0HB\\
 England.} 
\author{Stephen Doty}
\email{doty@math.luc.edu}
\address{Department of Mathematics and Statistics\\
Loyola University Chicago\\
6525 North Sheridan Road\\
Chicago\\ IL 60626\\ USA.}
 \author{Paul Martin} 
\email{ppmartin@maths.leeds.ac.uk}
\address{Department of Pure Mathematics\\
University of Leeds\\
Leeds\\
LS2 9JT\\
England.}
 \subjclass[2000]{Primary 20G05}
 \begin{abstract}We determine the blocks of the walled Brauer algebra
 in characteristic zero. These can be described in terms of orbits of
 the action of a Weyl group of type $A$ on a certain set of
 weights. In positive characteristic we give a linkage principle in
 terms of orbits of the corresponding affine Weyl group. We also
 classify the semisimple walled Brauer algebras in all
 characteristics.
 \end{abstract}

\maketitle

\section{Introduction}

The representation theories over $\CC$ of the symmetric group $\Sigma_r$
and the general linear group GL$_n(\CC)$ are related by Schur-Weyl
duality. This is the observation that the $r$th tensor product
$V^{\otimes r}$ of the natural representation $V$ of GL$_n(\CC)$ has
actions both by GL$_n(\CC)$ and by $\Sigma_r$ (the latter by place
permutation of the tensor factors) such that the image of each group
algebra under its action can be identified with the centraliser algebra
of the other \cite{weyl}.

The Brauer algebra $B_r(\delta)$ is an extension of $\CC\Sigma_r$
introduced to play the role of $\CC\Sigma_r$ in corresponding
dualities for symplectic and orthogonal groups \cite{brauer}. When
$\delta$ is a positive integer the duality is with O$_{\delta}(\CC)$;
when $\delta$ is a negative even integer the duality is with
Sp$_{-\delta}(\CC)$. However, the algebra $B_r(\delta)$ itself is
defined for all choices of $\delta\in\CC$.

The walled Brauer algebra (also known as the rational Brauer algebra)
$B_{r,s}(\delta)$ arises from a third version of Schur-Weyl
duality. For $\delta\in\NN$, consider the mixed tensor product
$V^{\otimes r}\otimes (V^*)^{\otimes s}$ of the natural representation
(and its dual) for GL$_{\delta}(\CC)$. There is a subalgebra
$B_{r,s}(\delta)$ of the Brauer algebra which acts on this product to
give a Schur-Weyl duality with the GL$_{\delta}(\CC)$ action. This
algebra was studied by Turaev \cite{turwall}, Koike \cite{koikewall},
and Benkart et al \cite{bchlls}. 

The symmetric group, Brauer, and walled Brauer algebras may be
considered over arbitrary fields. It is known that the respective
Schur-Weyl dualities continue to hold in types $A$ and $C$ (when the
field is infinite) and are expected to hold for other types. This is
well known in type $A$; see \cite{ddh} for type $C$.

Over $\CC$ the Brauer and walled Brauer algebras are isomorphic to the
corresponding centraliser algebra when $|\delta|>>0$, and hence must be
semisimple. However, for small values of $\delta$ the centraliser
algebra is only a quotient of the original algebra, and non-semisimple
cases can occur.

Until relatively recently the representation theory of the Brauer
algebra was ill understood. A precise semisimplicity criterion was
given by Rui \cite{ruibrauer} only very recently. However, in recent
work \cite{cdm,cdm2} it has been shown that the non-semisimple cases
have a rich combinatorial structure, which is controlled by the type
$D$ Weyl group; at present a structural explanation for this phenomena
is lacking.  The goal of the present paper is to analyse the
representation theory of the walled Brauer algebra in the same
manner. This will combine an application of the towers of recollement
formalism from \cite{cmpx}, generalised to the cellular setting, with
analogues of various explicit calculations for the Brauer algebra in
\cite{dhw}. Then we will introduce a geometric formulation of the
combinatorics obtained.

Section \ref{basics} introduces the walled Brauer algebra, and starts
to show how it is compatible with a cellular version of the basic
towers of recollement machinery. This is completed in Section
\ref{modules}, where the cell modules are introduced. 

In Section \ref{necsec} we give a necessary condition (for arbitrary
fields and parameter choices) for two simple modules to be in the same
block by considering the action of certain central elements in
$B_{r,s}$. This gives a necessary condition for semisimplicity; this
is shown to give a sufficient condition (except in obviously
non-semisimple cases) in Sections \ref{sufsec} and \ref{twosec} by
constructing certain homomorphisms between cell modules. Thus we
are able to give a complete semisimplicity criterion in Theorem
\ref{whenss}.

In Section \ref{blocksec} we are able to refine these results to give
a complete description of when two simple modules are in the same
block in characteristic zero.

For the Brauer algebra the combinatorial description of blocks in
\cite{cdm} has been reinterpreted \cite{cdm2} in terms of an action of
the Weyl group of type $D$ as in Lie theory. In Section \ref{geom1} we
will provide a similar reinterpretation for the walled Brauer, in
terms of the Weyl group of type $A$, but with an unusual choice of
dominant weights.

\begin{figure}[ht]
\includegraphics{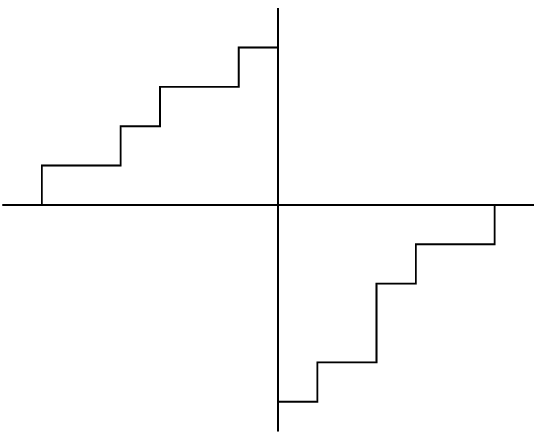}
\caption{A dominant weight for $B_{r,s}$}
\label{genform}
\end{figure}

Simple modules for the walled Brauer algebra $B_{r,s}$ can be labelled
by certain bipartitions $(\lambda^L,\lambda^R)$, and we choose to
identify these with elements of $\ZZ^{r+s}$ (with a fixed origin) by
embedding $\lambda^R$ normally but reflecting $\lambda^L$ about the
two axes passing through the origin. These dominant weights are thus
generalised partitions with Young diagram of the general form shown in
Figure \ref{genform}. As in Lie theory, the natural action of
$\Sigma_{r+s}$ on the set of weights has to be shifted, but in this
case the shift also depends on the parameter $\delta$.

In Section \ref{geom2} we give a necessary condition for two
weights to be in the same block in positive characteristic (i.e.~ a
linkage principle) by replacing $W$ by the corresponding affine Weyl
group. This is exactly analogous to the Brauer algebra case.

Although the combinatorics is a little more intricate (as the role of
the symmetric group for the Brauer algebra is replaced by a product of
symmetric groups), the actual proofs are rather simpler for the walled
Brauer algebra than for the Brauer algebra. This paper is also almost
entirely self-contained, requiring only a result of Halverson (on the
symmetric group content of walled Brauer modules) \cite{halvwall} from
the existing walled Brauer literature. Thus it may also be read as an
introduction to the methods used in the Brauer algebra papers
\cite{cdm} and \cite{cdm2}.

\section{The walled Brauer algebra}\label{basics}

Fix an algebraically closed field $k$ of characteristic $p\geq 0$, and
$\delta\in k$. For $r,s\in\NN$, the walled Brauer algebra
$B_{r,s}(\delta)$ can be defined as a subalgebra of the ordinary Brauer
algebra $B_{r+s}(\delta)$ in the following manner.

Recall that for $n\in\NN$ the \emph{Brauer algebra} $B_n(\delta)$ can
be defined in terms of a basis of partitions of
$\{1,\ldots,n,\overline{1},\ldots,\overline{n}\}$ into pairs. The
product $AB$ of two basis elements $A$ and $B$ is obtained by
representing each by a graph on $2n$ points, and identifying the
vertices $\overline{1},\ldots,\overline{n}$ of $A$ with the vertices
$1,\ldots,n$ of $B$ respectively. This produces a new graph on the
vertices $1,\ldots,n$ of $A$ and $\overline{1},\ldots,\overline{n}$ of
$B$, possibly together with some number ($t$ say) of connected
components not connected to any of these vertices. The product $AB$ is
then defined to be $\delta^tC$, where $C$ is the basis element
corresponding to the graph obtained by removing these connected
components.

It is usual to represent basis elements graphically by means of
diagrams with $n$ northern nodes numbered $1$ to $n$ from left to right, and
$n$ southern nodes numbered $\overline{1}$ to $\overline{n}$ from left
to right, where each node is connected to precisely one other by a
line. Edges connecting a northern and a southern node are called
\emph{propagating lines}, and the remainder are called \emph{northern}
or \emph{southern arcs}.

\begin{figure}[ht]
\includegraphics{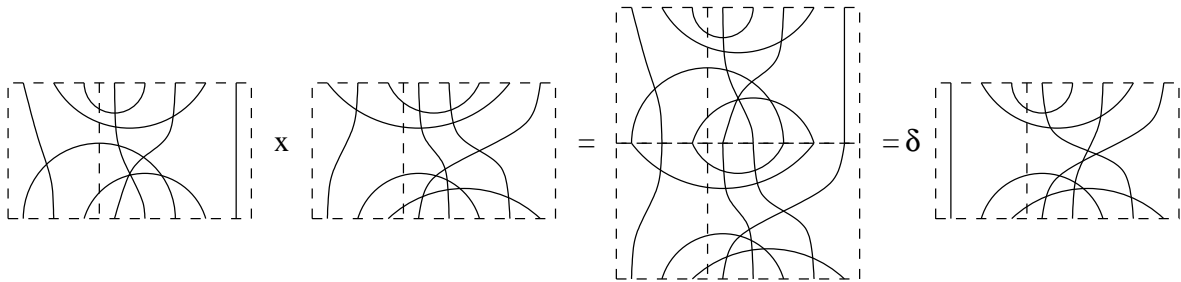}
\caption{The product of two basis elements in $B_{3,5}(\delta)$.}
\label{walldiagex}
\end{figure}

It is now easy to realise the \emph{walled Brauer algebra}
$B_{r,s}(\delta)$ as a subalgebra of the Brauer algebra
$B_{r+s}(\delta)$. Partition the basis diagrams with a wall separating
the first $r$ northern nodes and first $r$ southern nodes from the
remainder. Then the walled Brauer algebra is the subalgebra with basis
those Brauer diagrams such that no propagating edge crosses the wall,
and every northern or southern arc does cross the wall. It is easy to
verify that the space spanned by such diagrams is indeed a
subalgebra. Note that $B_{0,n}(\delta)\cong B_{n,0}(\delta)\cong
k\Sigma_n$, the group algebra of the symmetric group $\Sigma_n$ on $n$
letters.  An example of two walled Brauer diagrams and their product
is given in Figure \ref{walldiagex}.

We will show that the algebras $B_{r,s}(\delta)$ form cellular
analogues of the towers of recollement introduced
in\cite{cmpx}. Roughly, such a tower consists of a family of cellular
algebras related by inclusions and idempotent embeddings in a
compatible way, such that restriction and induction of cell modules is
well-behaved in the tower.  For general $k$ and $\delta$ we will
obtain a similar formalism will quasiheredity replaced by cellularity.

More precisely, there are six conditions labelled (A1--6) in
\cite{cmpx} which are required for a tower of recollement, and we will
consider each of these in turn.  All but (A2) (concerning
quasi-heredity) will turn out to hold (if suitably interpreted) for
arbitrary $k$ and $\delta$, and in the general case we will also be
able to replace (A2) by a cellular analogue.  Henceforth we will
suppress all $\delta$s in our notation when no ambiguity can occur.

\begin{figure}[ht]
\includegraphics{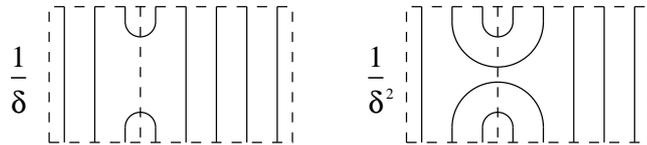}
\caption{The elements $e_{3,5}$ and $e_{3,5,2}$ in $B_{3,5}$.}
\label{idwall}
\end{figure}

Suppose that $k$ is arbitrary, with $r,s>0$ and $\delta\neq 0$, and
let $e_{r,s}\in B_{r,s}$ be $\delta^{-1}$ times the diagram with one
northern arc connecting $r$ and $r+1$, one southern arc connecting
$\overline{r}$ and $\overline{r+1}$ and all remaining edges being
propagating lines from $i$ to $\overline{i}$. (The element $e_{3,5}$
is illustrated in the left-hand side of Figure \ref{idwall}.) Clearly
$e_{r,s}$ is an idempotent in $B_{r,s}$.

\begin{figure}[ht]
\includegraphics{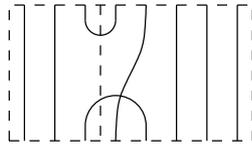}
\caption{The element $\tilde{e}_{3,5}$ in $B_{3,5}$.}
\label{altidem}
\end{figure}

If $\delta=0$ then we cannot define the idempotent $e_{r,s}$ as
above. However, if $r$ or $s$ is at least $2$ then we can define an
alternative idempotent $\tilde{e}_{rs}$ as illustrated in Figure
\ref{altidem}. 

Our first result is 

\begin{figure}[ht]
\includegraphics{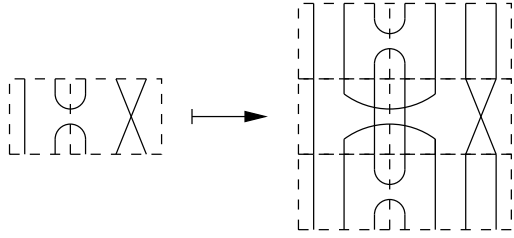}
\caption{An example of the action of the map $\Phi_{r,s}$.}
\label{eiso}
\end{figure}

\begin{prop}[A1]\label{A1}
If $\delta\neq 0$ then for each $r,s>0$ we have an algebra isomorphism 
$$\Phi_{r,s}:B_{r-1,s-1}\longrightarrow
e_{r,s}B_{r,s}e_{r,s}.$$
If $\delta=0$ and $r\geq 2$ or $s\geq 2$ we have an algebra
isomorphism 
$$\tilde{\Phi}_{r,s}:B_{r-1,s-1}\longrightarrow
\tilde{e}_{r,s}B_{r,s}\tilde{e}_{r,s}.$$
\end{prop}
\begin{proof}
We prove the first statement, the second is very similar.

Given a diagram $D$ in $B_{r-1,s-1}$ we define a new diagram $D'$ in
$B_{r,s}$ by adding two propagating lines immediately before and after
the wall in $D$, so that $r$ is connected to $\overline{r}$ and $r+1$
to $\overline{r+1}$. It is clear that the map taking $D$ to
$e_{r,s}D'e_{r,s}$ (as illustrated in Figure \ref{eiso}) is an
injective algebra homomorphism, and it is easy to verify that the
image is precisely $e_{r,s}B_{r,s}e_{r,s}$.
\end{proof}

\begin{rem}
The roles of $e_{r,s}$ and $\tilde{e}_{r,s}$ are very similar, and so
we will henceforth write $e_{r,s}$ for both types of idempotent (and
similarly write $\Phi$ for the isomorphisms in (A1)). This
will allow us to deal with the cases $\delta=0$ and $\delta\neq 0$
simultaneously. In proofs we will work with the original idempotent;
the obvious (trivial) modifications are left to the reader. 
\end{rem}

We wish to define a sequence of idempotents
$e_{r,s,i}$ in $B_{r,s}$. Set $e_{r,s,0}=1$, and for $1\leq i\leq
\min(r,s)$ set $e_{r,s,i}=\Phi_{r,s}(e_{r-1,s-1,i-1})$. Note that when
$\delta=0$ and $r=s$ the element $e_{r,r,r}$ is not defined.

To these elements we associate corresponding quotients
$B_{r,s,i}=B_{r,s}/B_{r,s}e_{r,s,i}B_{r,s}$. When $\delta\neq 0$ we
can give an alternative description of the $e_{r,s,i}$ (via
our explicit description of $\Phi_{r,s}$) as
$\delta^{-i}$ times the diagram with $i$ northern (respectively
southern) arcs connecting $r-t$ to $r+1+t$ (respectively
$\overline{r-t}$ to $\overline{r+1+t}$ for $0\leq t\leq i=1$ and the
remaining edges all propagating lines connecting $u$ to $\overline{u}$
for some $u$. The element $e_{3,5,2}$ is illustrated in the right-hand
side of Figure \ref{idwall}. A similar description can be given in the
case $\delta=0$.

We define the \emph{propagating vector} of a diagram $D$ to be the
pair $(a,b)$ where $D$ has $a$ propagating lines to the left of the
wall, and $b$ to the right.  Note that if we multiply
two diagrams with propagating vectors $(a_1,b_1)$ and $(a_2,b_2)$ then
the result must have propagating vector $(a,b)$ with $a\leq
\min(a_1,a_2)$ and $b\leq \min (b_1,b_2)$.

Set  $J_i= B_{r,s}e_{r,s,i}B_{r,s}$
and consider the sequence of ideals
\begin{equation}\label{qhchain}
\cdots\subset J_2\subset J_1\subset J_0=B_{r,s}.
\end{equation}

\begin{prop}\label{isec} The ideal $J_i$ has a basis of all diagrams
  with propagating vector $(a,b)$ for some $a\leq r-i$ and $b\leq
  s-i$. In particular the section
$$J_i/J_{i+1}$$
in the filtration (\ref{qhchain}) has a basis of all diagrams with
propagating vector $(r-i,s-i)$.
\end{prop}
\begin{proof}
This is a routine exercise (see for example \cite[Corollary
1.1]{msblob}).
\end{proof}

In particular  we have that 
\begin{equation}\label{topsec}
B_{r,s}/J_1\cong k(\Sigma_r\times\Sigma_s).
\end{equation}
We will denote $\Sigma_r\times\Sigma_s$ by $\Sigma_{r,s}$.

 We will need some basic facts about symmetric group
representations; details can be found in \cite{james}. For each
partition $\lambda$ of $n$, we can define a \emph{Specht module}
$S^{\lambda}$ for $k\Sigma_n$. We say that a partition
$\lambda=(\lambda_1,\ldots,\lambda_r)$ is \emph{$p$-regular}
if either $p>0$ and there is no $1\leq i\leq r$ such that 
$$\lambda_i=\lambda_{i+1}=\cdots=\lambda_{i+p}$$ or $p=0$.  Then the
heads $D^{\lambda}$ of the Specht modules $S^{\lambda}$ for $\lambda$
$p$-regular form a complete set of inequivalent simple
$k\Sigma_n$-modules.

As $k$ is algebraically closed (and so certainly a splitting field for
$\Sigma_r$ and $\Sigma_s$), the simple modules for $k\Sigma_{r,s}$ are
precisely those modules of the form $D\boxtimes D'$ where $D$ is a
simple $k\Sigma_r$-module, and $D'$ a simple $k\Sigma_s$-module
\cite[(10.33) Theorem]{cr1}, and
so can be labelled by pairs $(\lambda^L,\lambda^R)$ where $\lambda^L$
is a $p$-regular partition of $r$ and $\lambda^R$ is a $p$-regular
partition of $s$. We will denote the set of such pairs by $\Lambda^{r,s}_{reg}$.

If $p=0$ or $p>\max(r,s)$ then the group algebra
$k\Sigma_{r,s}$ is semisimple, and $\Lambda^{r,s}_{reg}$ consists of
all pairs of partitions of $r$ and $s$. If $(\lambda^L,\lambda^R)$ is such
a pair we will denote this by $(\lambda^L,\lambda^R)\vdash (r,s)$, and
denote the set of such by $\Lambda^{r,s}$.  We will call elements of 
$\Lambda^{r,s}$ \emph{weights}. We will say that $k$ is
\emph{$\Sigma_{r,s}$-semisimple} (or just \emph{$\Sigma$-semisimple}
whenever this does not cause confusion) when $p=0$ or $p>\max(r,s)$.

Let $\Lambda_{r,s}$ denote an indexing set for the simple
$B_{r,s}$-modules. From Proposition \ref{A1} (A1) we have an exact
\emph{localisation} functor
$$F_{r,s}:B_{r,s}\mbox{\rm -mod}\longrightarrow B_{r-1,s-1}\mbox{\rm
  -mod}$$ coming from the relevant idempotent, which takes a
$B_{r,s}$-module $M$ to $e_{r,s}M$. There is a
corresponding right exact \emph{globalisation} functor $G_{r-1,s-1}$
in the opposite direction which takes a $B_{r-1,s-1}$-module $N$ to
$B_{r,s}e_{r,s}\otimes_{e_{r,s}B_{r,s}e_{r,s}}N$. By standard
properties of localisation functors \cite{green} and (\ref{topsec}) we
have for $r,s>0$ that
$$
\Lambda_{r,s}=\Lambda_{r-1,s-1}\sqcup \Lambda^{r,s}_{reg}.
$$
As $B_{r,0}\cong B_{0,r}\cong k\Sigma_r$ we deduce

\begin{prop}\label{pairparts}
If $\delta\neq 0$ or
$r\neq s$, then
$$\Lambda_{r,s}=\coprod_{i=0}^{\min(r,s)}\Lambda^{r-i,s-i}_{reg}.$$
\end{prop}

It will be convenient to have the following alternative way of
describing walled Brauer diagrams, in terms of partial one-row
diagrams (confer for example \cite{gl}). Given a walled Brauer diagram
$D\in B_{r,s}$ with $t$ northern and $t$ southern arcs, we will write
$D=X_{v,w,\sigma}$ in the following manner. Let $v$ represent the
configuration of northern arcs in $D$, and $w$ represent the
configuration of southern arcs. Then $D$ is uniquely specified by
giving $\sigma\in\Sigma_{r-t,s-t}$ (regarded as a subset of
$\Sigma_{r+s}$ in the obvious way) such that $\sigma(i)=j$ if the
$i$th northern node on a propagating line is connected to node $j$. We
denote the set of elements $v$ arising thus by ${\mathcal V}_{r,s,t}$
(and by abuse of notation use the same set to refer to the elements
$w$ that arise), and call this the set of \emph{partial one-row
(r,s,t) diagrams}, or just \emph{partial one-row diagrams} when $r,s,$
and $t$ are clear from context.

The tower of recollement formalism in \cite{cmpx} is realised in the
context of quasihereditary algebras. However, it is easy to recast it
in the more general cellular algebra setting, albeit at the expense of
some additional hypotheses.

The notion of a cellular algebra was introduced by Graham and Lehrer
\cite{gl} in term of an involution and a basis with very special
properties. However, for our purposes the alternative (equivalent)
definition given later by K\"onig and Xi \cite{kxinfl} in terms of
ideals and iterated inflations (together with an involution) will ease
our exposition. (These two approaches also have a hybrid version in
the tabular framework \cite{margreentab}, although we will not
consider this here.) 

Given a $k$-algebra $C$, a $k$-vector space $V$, and a bilinear form
$\phi:V\otimes V\too C$, K\"onig and Xi define a (possibly nonunital)
algebra structure on $A_{C,V}^{\phi}=V\otimes V\otimes C$ by setting 
 the product of two 
basis elements to be
$$(a\otimes b\otimes x).(c\otimes d\otimes y)=a\otimes d\otimes
x\phi(b,c)y.$$ If $i$ is an involution on $C$ with
$i(\phi(v,w))=\phi(v,w)$ then there is an involution $j$ on
$A_{C,V}^{\phi}$ given by
$$j(a\otimes b\otimes x)=b\otimes a\otimes i(x).$$ The algebra
$A_{C,V}^{\phi}$ is called the \emph{inflation of $C$ along
$V$}. 

K\"onig and Xi also need to define algebra structures on sums of the
form $C\oplus D$ where $C$ is a (possibly nonunital) algebra and $D$
is a unital algebra, extending the two algebra structures and any
involutions which they possess in a compatible way. This is elementary
but rather involved; details can be found in
\cite{kxinfl, kxbrauer}. Iterating such constructions forms \emph{iterated
inflations}. The key result is that the inflation of a cellular
algebra is again cellular \cite[Proposition 3.3]{kxinfl}
In fact, carrying out this construction on full matrix algebras
gives precisely the class of cellular algebras \cite[Theorem
4.1]{kxinfl}.

In \cite[Section 5]{kxbrauer} these constructions were used to give a
simple proof that the Brauer algebra is cellular, by constructing it
as an iterated inflation of symmetric group algebras. We will modify
this argument to sketch a proof of a similar result for the walled
Brauer algebras involving the group algebras $k\Sigma_{m,n}$. This
result has been proved using tabular methods in \cite{margreentab}. An
explicit construction of a cellular basis can be found in \cite{enyang}.

There is an obvious involution $i$ on $B_{r,s}$ given by inverting
diagrams (so that northern nodes become southern nodes and vice
versa).

\begin{lem}\label{inflation}
For $l\geq 0$ the algebra $J_l/J_{l+1}$ is isomorphic to an inflation
$$V_l\otimes V_l\otimes k\Sigma_{r-l,s-l}$$ of $k\Sigma_{r-l,s-l}$
along a free $k$-module $V_l$ of rank $|{\mathcal V}_{r,s,l}|$, with
respect to some bilinear form $\phi$ (described below).
\end{lem}
\begin{proof} This is a very slight modification of the corresponding
  proof for Brauer algebras in \cite[Lemma 5.3]{kxbrauer}. Let $V_l$
  have basis ${\mathcal V}_{r,s,l}$, and let the map $$V_l\otimes
  V_l\otimes k\Sigma_{r-l,s-l}\too J_l/J_{l+1}$$ be given by $v\otimes
  w\otimes \sigma\longmapsto X_{v,w,\sigma}$. To define the value of
  $\phi(v,w)$ consider a product $X_{u,v,\sigma_1}X_{w,x,\sigma_2}$
  for some $u,x\in {\mathcal V}_{r,s,l}$ and
  $\sigma_1,\sigma_2\in\Sigma_{r-l,s-l}$. If this product does not
  have propagating vector $(r-l,s-l)$ then set
  $\phi(v,w)=0$. Otherwise $\phi(v,w)=\delta^t\sigma$ where $t$ is the
  number of closed loops in $X_{u,v,\sigma_1}X_{w,x,\sigma_2}$, and
  $\sigma$ is the unique permutation such that 
$$X_{u,v,\sigma_1}X_{w,x,\sigma_2}=\delta^tX_{u,x,\sigma_1\sigma\sigma_2}.$$
  Note that this definition is independent of the choice of
  $u,x,\sigma_1,\sigma_2$.  It is now easy to verify that we have the
  desired algebra isomorphism.
\end{proof}

Arguing exactly as in \cite{kxbrauer} one can then show that

\begin{prop}\label{Bisinf}
The walled Brauer algebra $B_{r,s}$ is an iterated inflation of 
group algebras of the form $\Sigma_{r-l,s-l}$ for $0\leq l\leq \min(r,s)$
along $V_l$.
\end{prop}

The group algebras of the symmetric groups are cellular \cite{gl}
(indeed they were the motivating example for cellularity), with cell
modules given by the Specht modules $S^{\lambda}$. From this follows

\begin{thm}\label{whencell}
(i) The walled Brauer algebra $B_{r,s}$ is cellular with a cell module
$\Delta_{r,s}(\lambda^L,\lambda^R)$ for each $(\lambda^L,\lambda^R)\in
\Lambda^{r-l,s-l}$ with $0\leq l\leq \min(r,s)$. \\ 
(ii) If
$\delta\neq 0$ or $r\neq s$ then the simple modules are indexed by all
pairs $(l,\lambda^L,\lambda^R)$ where $0\leq l\leq \min(r,s)$ and
$(\lambda^L,\lambda^R)\in\Lambda^{r-l,s-l}_{reg}$.\\ 
(iii) If $\delta=0$ and $r=s$ we get
the same indexing set for simples as in (ii), but with the single simple
corresponding to $l=\min(r,s)$ omitted.
\end{thm}
\begin{proof}
From the
basis definition in \cite{gl} (or see \cite[Proposition 6.15]{kxinfl})
it is clear that a cell basis for $k\Sigma_{r,s}$ can be obtained as a
product of cell bases for $k\Sigma_r$ and $k\Sigma_s$, and hence
$k\Sigma_{r,s}$ is cellular with cell modules of the form $M \boxtimes
N$ where $M, N$ are cell modules for $k\Sigma_r$, $k\Sigma_s$,
respectively. Part (i) now follows from Proposition \ref{Bisinf}.

Part (ii)  has already been
shown in Proposition \ref{pairparts}.  The modification needed when
$\delta=0$ and $r=s$ follows as in \cite[Corollary 5.8]{kxbrauer}, or
from the known structure of the cellular algebra $B_{1,1}$, as this is
identical to the Temperley-Lieb algebra TL$_2(0)$ .
\end{proof}

\begin{cor}[A2]\label{qhwhen}
If $k$ is $\Sigma$-semisimple, and either $\delta\neq 0$ or $\delta=0$
and $r\neq s$, then the algebra $B_{r,s}$ is quasi-hereditary, with
heredity chain induced by the idempotent $e_{r,s,i}$. In all other
cases $B_{r,s}$ is not quasihereditary.
\end{cor}
\begin{proof} This follows immediately from the fact that a cellular
  algebra is quasihereditary precisely when there are the same number
  of simples as cell modules. (The quasi-hereditary structure can also
  be proved directly as in \cite[Proposition 2.10]{cmpx}.)
\end{proof}

When $(\lambda^L,\lambda^R)\in\Lambda^{r-l,s-l}_{reg}$ for some $l\geq
0$ (with the same exception as in Theorem \ref{whencell}(iii)) we
shall denote the corresponding simple $B_{r,s}$-module by
$L_{r,s}(\lambda^L,\lambda^R)$. By standard cellular theory this 
arises as the head of the cell module $\Delta_{r,s}(\lambda^L,\lambda^R)$.

The tower of recollement formalism relies on the interplay between two
different ways of relating algebras: localisation/globalisation and
induction/restriction. Thus we also need a way of identifying one
walled Brauer algebra as a subalgebra of another. We will do this in
two different (but closely related) ways. 

For $r>0$ we may identify $B_{r-1,s}$ as a subalgebra of $B_{r,s}$,
and similarly for $B_{r,s-1}$ if $s>0$. There are a variety of ways of
doing this, but we will use

\begin{lem}[A3]
The map $\Psi_L$ (respectively $\Psi_R$) obtained by inserting a
propagating line immediately to the left (respectively right) of the
wall in a $B_{r-1,s}$ (respectively $B_{r,s-1}$) diagram extends to an
algebra inclusion of $B_{r-1,s}$ (respectively of $B_{r,s-1}$) inside
$B_{r,s}$.
\end{lem}

We thus have two restriction functors, $\res_{r,s}^L$ from
$B_{r,s}$-mod to $B_{r-1,s}$-mod and $\res_{r,s}^R$ from $B_{r,s}$-mod
to $B_{r,s-1}$-mod, and the corresponding right adjoint induction
functors $\ind_{r,s}^L$ from $B_{r,s}$-mod to $B_{r+1,s}$-mod and
$\ind_{r,s}^R$ from $B_{r,s}$-mod to $B_{r,s+1}$-mod. We will often
omit the subscripts from these functors when this is unambiguous.  Our
choice of algebra inclusions is motivated by the following
compatibility relation between restriction and localisation.

\begin{prop}[A4]\label{A4}
For all $k$ with $r,s>0$ (and $\delta\neq 0$ if $r=s=1$) we have that
$$B_{r,s}e_{r,s}\cong B_{r-1,s}$$ as a
$(B_{r-1,s},B_{r-1,s-1})$-bimodule, where the right action of
$B_{r-1,s-1}$ on $B_{r,s}e_{r,s}$ is given via the isomorphism in
(A1), and the left action of $B_{r-1,s}$ is given via the map
$\Psi_L$. There is a similar isomorphism
$$B_{r,s}e_{r,s}\cong B_{r,s-1}$$ as a
$(B_{r,s-1},B_{r-1,s-1})$-bimodule replacing $\Psi_L$ by $\Psi_R$.
\end{prop}
\begin{proof}
We will consider the first case, the second is similar.  Consider a
diagram D in $B_{r,s}e_{r,s}$. As a $(B_{r-1,s},B_{r-1,s-1})$-bimodule
this can be represented schematically as in the left-hand diagram in
Figure \ref{bimod}, where the shaded area above the diagram indicates
the northern nodes acted on via $B_{r-1,s}$ and below indicates the
southern nodes acted on via $B_{r-1,s-1}$. Note that the node $r$
(marked $A$) in the diagram is not acted on from above. We can convert
$D$ into a diagram for $B_{r-1,s}$ by removing the southern arc shown,
and deforming the line terminating at $A$ so that it terminates at the
point $B$ in the right-hand diagram in Figure \ref{bimod}. It is easy
to verify that this new diagram is indeed a walled Brauer diagram, and
lies in $B_{r-1,s}$. This gives the desired bimodule isomorphism.

\begin{figure}[ht]
\includegraphics{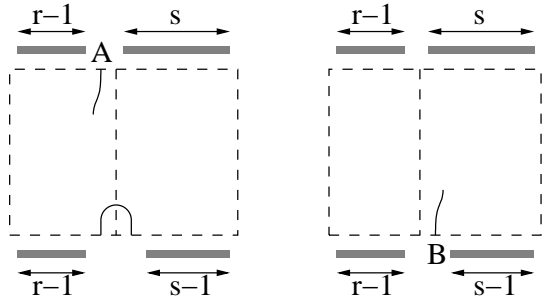}
\caption{Realising the bimodule isomorphism between $B_{r,s}e_{r,s}$
and $B_{r-1,s}$.}
\label{bimod}
\end{figure}
\end{proof}

Note that we have a choice of many different towers in our
construction, because at each stage we can chose either of the two
inclusions.

\section{Cell modules for the walled Brauer algebra}\label{modules}

To complete our verification of the tower of recollement axioms, and
their cellular analogues, we next analyse further the structure of the
cell modules.

We start by giving an explicit construction of the cell modules
(given a corresponding construction of Specht modules for the
symmetric groups).  For concreteness we will fix ${\mathcal
B}(\lambda^L,\lambda^R)$ to be the tensor product of the integral bases for
$S^{\lambda^L}$ and $S^{\lambda^R}$ given in \cite{jk}.

\begin{prop}\label{modbasis}
If $(\lambda^L,\lambda^R)\in\Lambda^{r-t,s-t}$ then the module 
$\Delta_{r,s}(\lambda^L,\lambda^R)$ has a basis given by
$$\{X_{v,1,id}\otimes x : v\in{\mathcal V}_{r,s,t}, \ x\in{\mathcal
B}(\lambda^L,\lambda^R)\}$$ where $X_{v,1,id}$ is a diagram with $r+s$
northern and  southern nodes, and $1$ denotes the (fixed)
southern half-diagram with arcs configured as in the southern half of
$e_{r,s,t}$.
\end{prop}
\begin{proof} This follows from the definition of cell ideals
  \cite[Definition 3.2 and the following proof]{kxcell} and the
  construction of cell ideals in inflations \cite[Section
  3.3]{kxinfl}, together with the explicit description of the
  inflations in Lemma \ref{inflation}.
\end{proof}

We will denote the space spanned by the elements $X_{v,1,id}$
with $v\in {\mathcal V}_{r,s,t}$ by $I_{r,s}^t$. 

Note that the action of a walled Brauer diagram on this basis is by
concatenation from above, except that products with too few
propagating lines are set equal to zero. As the action of the diagram
may induce a permutation $\sigma\in\Sigma_{r,s}$ of the propagating
lines this must be removed, and this is done by passing $\sigma$
through the tensor product to act on the basis element $x$ of
$S^{\lambda^L}\boxtimes S^{\lambda^R}$ in the natural manner. So if
$(\lambda^L,\lambda^R)\in\Lambda^{r-t,s-t}$ then 
$$\Delta_{r,s}(\lambda^L,\lambda^R)\cong
I_{r,s}^t\otimes_{\Sigma_{r-t,s-t}}(S^{\lambda^L}\boxtimes
S^{\lambda^R}).$$

\begin{cor}\label{fix}
For $(\lambda^L,\lambda^R)\in\Lambda^{r-t,s-t}$ the cell module
$\Delta_{r,s}(\lambda^L,\lambda^R)$ can be identified with the module
$$B_{r,s}e_{r,s,t}\otimes_{e_{r,s,t}B_{r,s}e_{r,s,t}}S^{\lambda^L}\boxtimes
S^{\lambda^R}$$ (when $e_{r,s,t}$ exists).
\end{cor}
\begin{proof}
This follows as in the proof of \cite[Proposition 2.10]{cmpx}.
\end{proof}

 We say that $(\lambda^L,\lambda^R)\leq (\mu^L,\mu^R)$ if
$(\lambda^L,\lambda^R)= (\mu^L,\mu^R)$ or
$(\lambda^L,\lambda^R)\in\Lambda^{a,b}$ and
$(\mu^L,\mu^R)\in\Lambda^{a-t,b-t}$ for some $0\leq t<\min(a,b)$. As
this ordering is compatible with the cellular structure, all
composition factors of $\Delta_{r,s}(\lambda^L,\lambda^R)$ are
labelled by weights $(\mu^L,\mu^R)$ with $(\mu^L,\mu^R)\leq
(\lambda^L,\lambda^R)$. It follows from the construction that
$$\Delta_{r,s}(\lambda^L,\lambda^R)\cong S^{\lambda^L}\boxtimes
S^{\lambda^R}$$ if $(\lambda^L,\lambda^R)\in \Lambda^{r,s}$, the lift
of a Specht module for the quotient algebra
$$B_{r,s}/J_1\cong k\Sigma_{r,s}.$$ 

As our globalisation and localisation functors (when they exist) are
compatible with the cell chain, we have by Corollary \ref{fix}
 that
\begin{equation}\label{Gident}
G_{r,s}(\Delta_{r,s}(\lambda^L,\lambda^R))\cong
\Delta_{r+1,s+1}(\lambda^L,\lambda^R)
\end{equation}
for all $(\lambda^L,\lambda^R)\in\Lambda_{r,s}$
and
$$F_{r,s}(\Delta_{r,s}(\lambda^L,\lambda^R))\cong \left\{\begin{array}{ll}
\Delta_{r-1,s-1}(\lambda^L,\lambda^R)&\wif (\lambda^L,\lambda^R)
\in\Lambda_{r-1,s-1}\\
0& \wif (\lambda^L,\lambda^R)\in\Lambda^{r,s}.\end{array}\right.
$$ 
As $F_{r,s}$ is exact we also have that
$$F_{r,s}(L_{r,s}(\lambda^L,\lambda^R))\cong \left\{\begin{array}{ll}
L_{r-1,s-1}(\lambda^L,\lambda^R)&\wif (\lambda^L,\lambda^R)
\in\Lambda_{r-1,s-1}\\
0& \wif (\lambda^L,\lambda^R)\in\Lambda^{r,s}.\end{array}\right.
$$ 

The compatibility of induction/restriction with
localisation/globalisation given in Proposition \ref{A4} (A4)
immediately implies that
\begin{equation}\label{resG}
\res^{\dagger}(G_{r,s}(\Delta_{r,s}(\lambda^L,\lambda^R))\cong
\ind^{\ddag}\Delta_{r,s}(\lambda^L,\lambda^R)
\end{equation} for all
$(\lambda^L,\lambda^R)\in\Lambda_{r,s}$ where $(\dagger,\ddag)$ represents
either $(L,R)$ or $(R,L)$.

Note that the only case where localisation and globalisation functors
do not exist occurs when $\delta=0$ and $r=s$. In this case we do not
have $F_{1,1}$ and $G_{0,0}$.

The remaining two axioms (A5) and (A6) concern the behaviour of
cell modules under induction and restriction. 

We first consider the restriction rules for cell modules. For this
we need to consider the action of certain special elements in the
walled Brauer algebra.  Write $E_{i,j}$ for the walled Brauer diagram
with edges between $t$ and $\overline{t}$ for $t\neq i,j$, and arcs
between $i$ and $j$, and $\overline{i}$ and $\overline{j}$. Note that
$B_{r,s}$ is generated by the elements $E_{i,j}$ (with $1\leq i\leq r$
and $r+1\leq j\leq r+s$) and the group $\Sigma_{r,s}$ (identified with
the set of diagrams with no northern or southern arcs). (In fact,
$B_{r,s}$ is generated by $\Sigma_{r,s}$ together with just 
one $E_{i,j}$.)

Consider the action of $E_{i,j}$ on an element $X_{w,1,id}\otimes x\in
\Delta_{r,s}(\lambda^L,\lambda^R)$. There are four possible cases:

(a) $i$ and $j$ are connected in $w$. In this case the action of
$E_{i,j}$ creates a closed loop while leaving the underlying diagram
unchanged. Hence
\begin{equation}\label{acase}
E_{i,j}(X_{w,1,id}\otimes x)=\delta(X_{w,1,id}\otimes x).
\end{equation}

(b) $i$ and $j$ are free vertices in $w$. In this case the action of
$E_{i,j}$ creates an extra northern arc, hence reducing the number of
propagating lines by two. Hence
\begin{equation}\label{bcase}
E_{i,j}(X_{w,1,id}\otimes x)=0.
\end{equation}

\begin{figure}[ht]
\includegraphics{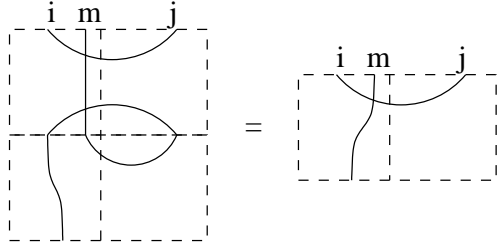}
\caption{A schematic example of the action of $E_{i,j}$ in case (c)}
\label{casec}
\end{figure}

(c) One of the vertices $i$ and $j$ is free, and the other is joined
to some vertex $m$ in $w$. Suppose that $i$ is the free vertex. Then
the action of $E_{i,j}$ is illustrated schematically in Figure
\ref{casec} (where we have omitted all lines which do not concern
us). From this it is clear that
\begin{equation}\label{ccase}
E_{i,j}(X_{w,1,id}\otimes x)=(i,m)(X_{w,1,id}\otimes x)
\end{equation} 
where $(i,m)$ denotes the transposition swapping $i$ and $m$ in
$\Sigma_{r,s}\subset B_{r,s}$. A similar result holds if we
reverse the roles of $i$ and $j$.

\begin{figure}[ht]
\includegraphics{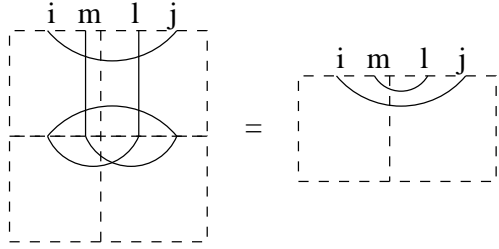}
\caption{A schematic example of the action of $E_{i,j}$ in case (d)}
\label{cased}
\end{figure}

(d) $i$ and $j$ are unconnected, but not free vertices in $w$. Say $i$
is joined to $l$ and $j$ is joined to $m$. This case is illustrated
schematically in Figure \ref{cased}. From this it is clear that
\begin{equation}\label{dcase}E_{i,j}(X_{w,1,id}\otimes x)
=(i,m)(X_{w,1,id}\otimes x)
=(j,l)(X_{w,1,id}\otimes x).\end{equation}

We will denote the Young diagram associated to a partition $\lambda$
by $[\lambda]$. For a partition $\lambda$, recall that the set of
\emph{removable boxes} are those which can be removed (singly) from
$[\lambda]$ such that the result is the Young diagram of a
partition. Similarly the set of \emph{addable boxes} are those which
can be added (singly) to $[\lambda]$ such that the result is the Young
diagram of a partition. We denote these sets by $\remo(\lambda)$ and
$\add(\lambda)$ respectively, and given a box in $\remo(\lambda)$ or
$\add(\lambda)$ denote the associated partition obtained by addition
or subtraction by $\lambda\pm\square$.

Given two partitions $\lambda$ and $\mu$ of $n$, we say that $\mu$ is
\emph{dominated} by $\lambda$ (written $\mu\lhd\lambda$) if for all
$i\geq 1$ we have
$$\sum_{j=1}^i\mu_j\leq\sum_{j=1}^i\lambda_j.$$ That is, the Young
diagram for $\lambda$ can be obtained from that for $\mu$ by moving
some of the boxes to earlier rows in the diagram. We will extend this
to give a partial order on pairs of partitions in a very
restricted form by saying that $(\lambda^L,\lambda^R)\lhd^L(\mu^L,\mu^R)$ if
$\lambda^R=\mu^R$ and $\lambda^L\lhd \mu^L$, and similarly for
$\lhd^R$ reversing the roles of $L$ and $R$.

Given a family of modules $M_i$ we will write $\biguplus_i M_i$ to
 denote some module with a filtration whose quotients are exactly the
 $M_i$, each with multiplicity $1$. This is not uniquely defined as a
 module, but the existence of a module with such a filtration will be
 sufficient for our purposes.

 We can now prove 

\begin{thm}\label{resrule}
(i) Suppose that $(\lambda^L,\lambda^R)\in\Lambda^{r-t,s-t}$. If $t=0$ then
$$\res_{r,s}^L\Delta_{r,s}(\lambda^L,\lambda^R)\cong
\biguplus_{\square\in\remo(\lambda^L)}
\Delta_{r-1,s}(\lambda^L-\square,\lambda^R).$$
If $t>0$ then we have a short exact sequence
$$0\too\!\!  \biguplus_{\square\in\remo(\lambda^L)}\!
\Delta_{r-1,s}(\lambda^L-\square,\lambda^R) \too
\res_{r,s}^L\Delta_{r,s}(\lambda^L,\lambda^R) \too\!\!
\biguplus_{\square\in\add(\lambda^R)}\!
\Delta_{r-1,s}(\lambda^L,\lambda^R+\square) \too 0$$ where the
left-hand sum equals $0$ if $\lambda^L=\emptyset$. \\
(ii) In each of the
filtered modules which arise in (i), the filtration can be chosen so that
the weights labelling successive quotients are ordered by $\lhd^L$ or
$\lhd^R$, with the top quotient
maximal among these. When $k$ is
$\Sigma$-semisimple the $\uplus$ all become direct sums.\\
(iii) There is a similar result for $\res_{r,s}^R$ replacing
$\remo(\lambda^L)$  by $\remo(\lambda^R)$  
and $\add(\lambda^R)$  by $\add(\lambda^L)$. 
\end{thm}
\begin{proof} We prove the result for $\res_{r,s}^L$; the right-hand case is 
similar. Our proof is very similar to that for the ordinary Brauer
algebra in \cite[Theorem 4.1]{dhw}.

Let $W$ be the subspace of $\Delta_{r,s}(\lambda^L,\lambda^R)$ spanned
by elements of the form $X_{w,1,id}\otimes x$ where the node in
$X_{w,1,id}$ numbered $r$ is on a propagating line. Recall our
realisation of $B_{r-1,s}$ inside $B_{r,s,}$. It is clear that the
elements of $\Sigma_{r-1,s}$, and the $E_{i,j}$ with $1\leq i<r$ and
$r+1\leq j\leq r+s$ preserve the space $W$, and hence $W$ is a
$B_{r-1,s}$-submodule.

We will show that 
$$W\cong
\biguplus_{\square\in\remo(\lambda^L)}
\Delta_{r-1,s}(\lambda^L-\square,\lambda^R).$$
When $t=0$ the space $W$ is the whole of $\Delta_{r,s}(\lambda^L,\lambda^R)$,
and so this will complete the proof of (i) in that case.  We have
\begin{multline}\label{resay}
I_{r-1,s}^t\otimes_{\Sigma_{r-1-t,s-t}}
\res_{\Sigma_{r-t-1,s-t}}^{\Sigma_{r-t,s-t}}
(S^{\lambda^L}\boxtimes
S^{\lambda^R})\\\cong 
I_{r-1,s}^t\otimes_{\Sigma_{r-1-t,s-t}}\biguplus_{\square\in\remo(\lambda^L)}
S^{\lambda^L-\square} \boxtimes S^{\lambda^R}\\
\cong \biguplus_{\square\in\remo(\lambda^L)}I_{r-1,s}^t 
\otimes_{\Sigma_{r-1-t,s-t}} S^{\lambda^L-\square} \boxtimes
S^{\lambda^R}
\cong
\biguplus_{\square\in\remo(\lambda^L)}
\Delta_{r-1,s}(\lambda^L-\square,\lambda^R).
\end{multline}
By \cite[Theorem 9.3]{james}, the restriction of a Specht module
satisfies condition (ii) above, which will thus be inherited by our
filtered module.
Thus it is enough to  show that $W$ is isomorphic to 
$I_{r-1,s}^t\otimes_{\Sigma_{r-1-t,s-t}}
\res_{\Sigma_{r-t-1,s-t}}^{\Sigma_{r-t,s-t}}
(S^{\lambda^L}\boxtimes
S^{\lambda^R})$.

Given a diagram $d=X_{w,1,id}$ with a propagating line from node
$r$, let $\bar{\phi}(d)$ be the diagram obtained by deleting this
line. We claim that the map
$$\phi:W\longrightarrow
I_{r-1,s}^t\otimes_{\Sigma_{r-1-t,s-t}}
\res_{\Sigma_{r-t-1,s-t}}^{\Sigma_{r-t,s-t}}
(S^{\lambda^L}\boxtimes
S^{\lambda^R})$$ given by $d\otimes x\longmapsto \bar{\phi}(d)\otimes x$
provides the desired isomorphism. Note that as vector spaces the
isomorphism is clear, as $\bar{\phi}$ is a bijection. Thus it is
enough to show that $\phi$ commutes with the action of
$\Sigma_{r-1,s}$ and $E_{i,j}$ with $1\leq i<r$ and $r+1\leq j\leq
r+s$. That $\phi$ commutes with the $\Sigma_{r-1,s}$ action is clear,
and by our discussion of the cases (a)--(d) above, this is also the
case for the $E_{i,j}$ action.

It remains to show that when $t>0$ the quotient 
$$V=\Delta_{r,s}(\lambda^L,\lambda^R)/W\cong
\biguplus_{\square\in\add(\lambda^R)}
\Delta_{r-1,s}(\lambda^L,\lambda^R+\square).$$
Arguing as in (\ref{resay}) (using \cite[17.14]{james} which gives
(ii) for induction of Specht modules),
it is easy to see that it is enough to
show that
$$V\cong I_{r-1,s}^{t-1}\otimes_{\Sigma_{r--t,s-t+1}}
\ind_{\Sigma{r-t,s-t}}^{\Sigma_{r-t,s-t+1}}(S^{\lambda^L}\boxtimes S^{\lambda^R}).
$$

We will need the following explicit realisation of the induced
module. Let $a=r-t$ and $b=s-t$.  For $a+1\leq i\leq a+b$ let $\tau_i$
be the transposition $(i,a+b+1)$, and let $\tau_{a+b+1}=1$. These
elements form a set of coset representatives for
$\Sigma_{a,b+1}/\Sigma_{a,b}$. If $A$ is a basis for
$S^{\lambda^L}\boxtimes S^{\lambda^R}$ then $A\times \{a+1,\ldots,a+b+1\}$ can
be regarded as a basis for
$\ind_{\Sigma{a,b}}^{\Sigma_{a,b+1}}(S^{\lambda^L}\boxtimes
S^{\lambda^R})$. The action of $\sigma\in\Sigma_{a,b+1}$ on such a basis
element $(v,j)$ is given by
\begin{equation}\label{indaction}
\sigma(v,j)=\left((\tau_l\sigma\tau_j)v,l\right)
\end{equation} where $l$ is the
unique value such that $\tau_l\sigma\tau_j\in\Sigma_{a,b}$.

Consider the group algebra $k\Sigma_{a,b+1}$ as a subset of the walled
Brauer algebra $B_{a,b+1}$ in the usual way, and suppose that we use
a diagrammatic notation for representing the action of $\Sigma_{a,b}$ on
$S^{\lambda^L}\boxtimes S^{\lambda^R}$ where the action of the group
$\Sigma_{a,b}$ is via the $a+b$ propagating lines. We will wish
to represent the action of $\Sigma_{a,b+1}$ on $\ind
S^{\lambda^L}\boxtimes S^{\lambda^R}$ in a similar manner, with the
aid of a dummy node at the right-hand end of the southern edge of a
diagram. Given a basis element $v$ for $S^{\lambda^L}\boxtimes
S^{\lambda^R}$, denote the element $(v,j)$ by a diagram of the form
shown in Figure \ref{realise}, where node $*$ denotes the $j$th
northern node (from the left) on a propagating line, all northern arcs
have been suppressed for simplicity, and the only lines that cross
occur when a line crosses the line from $*$. The basis element of $v$
is symbolically attached to the first $a+b$ southern nodes on
propagating lines as usual, and the final node is a dummy node. It is
now a routine exercise to verify that this representation coincides
with the action in (\ref{indaction}) under concatenation of diagrams,
and so is a diagrammatic realisation of the induced module.

\begin{figure}[ht]
\includegraphics{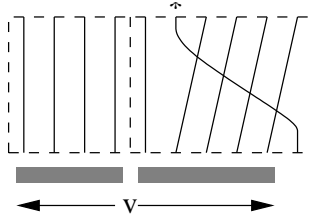}
\caption{Realising the element $(v,j)$ as a diagram.}
\label{realise}
\end{figure}

Given a diagram $d=X_{w,1,id}$ where $w$ has a northern arc connecting
node $r$ to some other node $i$ say, we form a new diagram $\psi(d)$
by deleting this arc, and replacing it by a propagating line from $i$
to a new node at the right-hand end of the southern edge. Note that
this new edge may cross over some of the original propagating lines,
and that there is a unique permutation in $\Sigma_{a,b+1}$ which
transforms it into an element of the form $X_{w,1,id}$. We
claim that the map
$$\psi:V\longrightarrow 
I_{r-1,s}^{t-1}\otimes_{\Sigma_{r-t,s-t+1}}
\ind_{\Sigma{a,b}}^{\Sigma_{a,b+1}}(S^{\lambda^L}\boxtimes S^{\lambda^R})$$
given by
$$d\otimes x\longmapsto \psi(d)\otimes(x,a+b+1)$$ (where we represent
$(x,a+b+1)$ diagrammatically as above) gives the desired
$B_{r-1,s}$-isomorphism.  First note that as $\bar{\psi}$ is a $(b+1)$
to $1$ map, and
$$\dim \ind_{\Sigma{a,b}}^{\Sigma_{a,b+1}}(S^{\lambda^L}\boxtimes
S^{\lambda^R})=(b+1)\dim (S^{\lambda^L}\boxtimes
S^{\lambda^R})$$
the map $\psi$ is between two spaces of the same dimension.

We next show that $\psi$ is onto. A basis for $I_{r-1,s}^{t-1}\otimes
\ind_{\Sigma{a,b}}^{\Sigma_{a,b+1}}(S^{\lambda^L}\boxtimes
S^{\lambda^R})$ is given by $\{X_{w,1,id}\otimes (x,j)\}$ where $x$
runs over a basis for $S^{\lambda^L}\boxtimes S^{\lambda^R}$. From the
explicit description of our action on the induced module, we can
easily see that there exists some second basis element $x'$ for
$S^{\lambda^L}\boxtimes S^{\lambda^R}$ such that $\sigma_j(x,j)=(x',a+b+1)$,
as $\tau_{a+b+1}\sigma_j\tau_j\in\Sigma_{a,b}$. Therefore
$$X_{w,1,id}\otimes (x,j)=X_{w,1,id}\otimes \sigma_j(x',a+b+1)=
X_{w,1,\rho}\otimes (x',a+b+1).$$ But it is easy to construct an element
$d$ such that $\psi(d)=X_{w,1,\rho}$, and hence $X_{w,1,id}\otimes
(x,j)$ is in the image of $\psi$.

Thus we will be done once we have shown that $\psi$ is a
$B_{r-1,s}$-homomorphism. By considering the diagrammatic realisation
of the induced module, it is easy to see that this is a
$\Sigma_{r-1,s}$-homomorphism, so it will be enough to show that the
action of the $E_{i,j}$ commutes with $\psi$.

Consider a basis element $b\otimes x$ in $V$, and an element $E_{i,j}$
with $1\leq i<r$ and $r+1\leq j\leq r+s$. The northern node $r$ in $b$
is connected to some northern node $t$ with $r+1\leq t\leq r+s$. If
$j\neq t$ then it is clear that the action of $E_{i,j}$ commutes with
$\psi$. If $j=t$ and the northern node $i$ in $b$ is on a northern arc
then the result is also clear. 

Finally, suppose that $j=t$ and the northern node $i$ in $b$ is on a
propagating line. Then in $E_{i,j}(b\otimes x)$ the northern node $r$
is on a propagating line, and hence this product is in $W$ (and so is
$0$ in $V$). The action of $E_{i,j}$ on $\psi(d)\otimes (x,a+b+1)$ is
clearly $0$ as the diagram obtained has an extra arc, and hence too
few propagating lines. Thus $E_{i,j}$ commutes with the action of
$\psi$ in all cases, and so we are done. 
\end{proof}

\begin{rem}\label{better}
In the non-$\Sigma$-semisimple cases, we get partial direct sum
decompositions for the filtered modules in Theorem \ref{resrule}
associated to the corresponding block decompositions for the symmetric
group. This follows from (\ref{resay}) and its analogue for induction.
\end{rem}

If we are in the quasi-hereditary case then we have that
$\Lambda_{r,s}=\coprod_{t=0}^{\min(r,s)}\Lambda^{r-t,s-t}$.  We will
write $\Lambda_{r,s}^m$ for the set $\Lambda^{r-t,s-t}$ with
$m=r+s-2t$ when regarded as a subset of $\Lambda_{r,s}$ in this
way. (Note that $m=|\mu^L|+|\mu^R|$ for $(\mu^L,\mu^R)\in
\Lambda^{r-t,s-t}$.)  Given a $B_{r,s}$-module $M$ with a filtration
by cell modules (which we call a \emph{cell filtration}), we would
like to be able to define the \emph{support} $\supp_{r,s}(M)$ of $M$
to be the set of labels of cell modules in such a filtration. In the
quasi-hereditary case this is well-defined as standard modules form a
basis for the Grothendieck group, however for general cellular
algebras it may depend on the filtration chosen. However, for our
purposes we will only ever need that there exists some filtration with
a certain support.

We start with a weaker version of Theorem \ref{resrule}, verifying the
fifth tower axiom:

\begin{cor}[A5]
For each $(\lambda^L,\lambda^R)\in\Lambda_{r,s}^m$ we have that
$\res^L(\Delta_{r,s}(\lambda^L,\lambda^R))$ has a cell filtration,
such that
$$\supp(\res^L(\Delta_{r,s}(\lambda^L,\lambda^R)))
\subseteq \Lambda_{r-1,s}^{m-1}\sqcup\Lambda_{r-1,s}^{m+1}$$
and a similar result for $\res^R$.
\end{cor}

From (A4) we have that 
$$\res^{\dagger}G(\Delta_{r,s}(\lambda^L,\lambda^R))\cong
\ind^{\ddag}\Delta_{r,s}(\lambda^L,\lambda^R)$$
where $(\dagger,\ddag)$ equals $(L,R)$ or $(R,L)$. Combining this with (A5) and
(\ref{Gident}) we obtain the analogue of Theorem \ref{resrule} for
induction:

\begin{cor}\label{indrule}
Suppose that $(\lambda^L,\lambda^R)\in\Lambda^{r-t,s-t}$.  Then we
have a short exact sequence
$$0\too\!\!
\biguplus_{\square\in\remo(\lambda^L)}\!
\Delta_{r,s+1}(\lambda^L-\square,\lambda^R)
\too \ind_{r,s}^R\Delta_{r,s}(\lambda^L,\lambda^R) \too\!\!
\biguplus_{\square\in\add(\lambda^R)}\!
\Delta_{r,s+1}(\lambda^L,\lambda^R+\square)
\too 0$$
where the left-hand sum equals $0$ if $\lambda^L=\emptyset$. 
Each of the filtered modules arising has filtration compatible with
the $\lhd^L$ or $\lhd^R$ order as in Theorem \ref{resrule}, and when $k$ is
$\Sigma$-semisimple the $\uplus$ all become direct sums.
 There is a similar result for $\ind_{r,s}^L$ replacing
$\remo(\lambda^L)$  by $\remo(\lambda^R)$  
and $\add(\lambda^R)$  by $\add(\lambda^L)$.
\end{cor}

Note that the roles of $R$ and $L$ are reversed for induction and
restriction rules. We also have 

\begin{cor}
For each $(\lambda^L,\lambda^R)\in\Lambda_{r,s}^m$ we have that
$\ind^R(\Delta_{r,s}(\lambda^L,\lambda^R))$ has a cell filtration,
such that 
$$\supp(\ind^R(\Delta_{r,s}(\lambda^L,\lambda^R)))
\subseteq \Lambda_{r,s+1}^{m-1}\sqcup\Lambda_{r,s+1}^{m+1}$$
and a similar result for $\ind^L$.
\end{cor}

Finally note that it is clear from the precise form of the induction
rules that

\begin{cor}[A6] If $r>0$ then for
 each $(\lambda^L,\lambda^R)\in\Lambda_{r,s}^{r+s}$ there exists
 $(\mu^L,\mu^R)\in\Lambda_{r-1,s}^{r+s-1}$, and a cell filtration of
 $\ind_{r-1,s}^L\Delta_{r-1,s}(\mu^L,\mu^R)$, such that 
$$(\lambda^L,\lambda^R)\in
\supp(\ind_{r-1,s}^L\Delta_{r-1,s}(\mu^L,\mu^R)).$$
There is a similar result if $s>0$ involving $\ind_{r,s-1}^R$.
\end{cor}

The tower of recollement machinery was introduced in \cite{cmpx} in
part as an organisational tool to reduce the analysis of
representations of towers of
algebras to certain special cases. To this extent the adaptations
above with cellularity instead of quasi-heredity are sufficient.

We would like to have a tower of recollement (or cellular analogue)
containing $B_{r,s}$ in the sense of \cite{cmpx}. As stated there the
tower depends on one indexing label rather than two, and so we will
need to make choices for our algebra inclusions in a consistent
way. Suppose that $r\geq s$ (the case $r<s$ is similar). Let $t=r-s$,
and set $A_0=B_{t,0}$. For $n>0$ set
$$A_n=\left\{\begin{array}{ll} B_{t+u,u+1} & \wif n=2u+1\\ B_{t+u,u} &
\wif n=2u.\end{array}\right.$$ The functors $F$ and $G$ corresponding
to each choice of $A_n$ are compatible with this choice of algebras,
and give functors $F_n$ from $A_n$-mod to $A_{n-2}$-mod and $G_n$ from
$A_n$-mod to $A_{n+2}$-mod. Modules in $A_n$-mod will be identified by
the subscript $n$ instead of the corresponding pair $r,s$. Choosing
$\ind_n$ and $\res_n$ alternately left and right as $n$ increases (so
that they go from $A_{n-1}$-mod to $A_n$-mod and vice versa) it is
easy to check that in the quasi-hereditary case (A1--6) now follow
exactly as in \cite{cmpx}, with $\Lambda_n=\Lambda_{a,b}$ and
$\Lambda_n^m=\Lambda_{a,b}^{m-t}$ if $A_m=B_{a,b}$.

We have the following cellular version of \cite[Theorem 3.7]{cmpx} for
$B_{r,s}$.

\begin{thm}\label{reduceto} 
(i) For all $(\lambda^L,\lambda^R)\in\Lambda_n^m$ and
$(\mu^L,\mu^R)\in\Lambda_n^l$ we have
$$\Hom(\Delta_n(\lambda^L,\lambda^R),\Delta_n(\mu^L,\mu^R))\cong
\left\{\begin{array}{ll}
\Hom(\Delta_m(\lambda^L,\lambda^R),\Delta_m(\mu^L,\mu^R)) & \wif l\leq
m\\ 0 & \otherwise.\end{array}\right.$$ (ii) The algebras $A_n$ are
semisimple for $0\leq n\leq N$ if and only if they are
quasi-hereditary and for all $0\leq n\leq N$ and pairs of weights
$(\lambda^L,\lambda^R)\in\Lambda_{n}^{n}$ and
$(\mu^L,\mu^R)\in\Lambda_{n}^{n-2}$ we have
$$\Hom(\Delta_n(\lambda^L,\lambda^R),\Delta_n(\mu^L,\mu^R))\cong 0.$$
\end{thm}
\begin{proof}
Part (i) follows exactly as in the proof of \cite[Theorem
  1.1]{cmpx}. For part (ii), note that if a cellular algebra is not
quasi-hereditary then it cannot be semisimple \cite[(3.8)
  Theorem]{gl}, so it is enough to consider the quasi-hereditary
case. This has already been proved in \cite[Theorem 1.1]{cmpx}.
\end{proof}

\section{A necessary condition for blocks}\label{necsec}

The principal aim of this section is to give a necessary condition for
two weights to label simple $B_{r,s}$ modules in the same block. (We
will abuse terminology and say that the weights themselves are in the
same block.)  This is closely modelled on a similar result for the
Brauer algebra in \cite[Theorem 3.3]{dhw} (as interpreted in
\cite[Proposition 4.2]{cdm}). Throughout this section $k$ and $\delta$
are arbitrary.

Let $\lambda$ be a partition. For a box $d$ in row $i$ and column $j$
of the Young diagram $[\lambda]$ we set $c(d)=j-i$, the \emph{content}
of $d$. We set
$$T_{r,s}=\sum_{\substack{
1\leq i\leq r,\\ r+1\leq j\leq r+s}}E_{i,j}.$$

\begin{lem}\label{keyscalar}
Let $(\lambda^L,\lambda^R)\in\Lambda^{r-t,s-t}$. For all
$y\in\Delta_{r,s}(\lambda^L,\lambda^R)$ we have that
$$T_{r,s}y= \Bigl(
t\delta-\sum_{d\in[\lambda^L]}c(d)-\sum_{d\in[\lambda^R]}c(d)+
\sum_{1\leq i<m\leq r}(i,m)+
\sum_{r< i<m\leq r+s}(i,m)\Bigr)y$$ where $(i,m)$ denotes the
element of $\Sigma_{r,s}$ which transposes $i$ and $m$.
\end{lem}
\begin{proof}
First suppose that $k=\CC$.  It is enough to consider the case where
$y=X_{w,1,id}\otimes x$ where $w\in {\mathcal V}_{r,s,t}$ and $x\in
S^{\lambda^L}\boxtimes S^{\lambda^R}$. Recall that the action of
$E_{i,j}$ on such an element is given by one of the four cases
(a)--(d) given by equations (\ref{acase}--\ref{dcase}).

We need to determine the contribution of each of the four cases to our
final sum. As there are $t$ northern arcs, there are precisely $t$
distinct pairs $(i,j)$ in case (a), and so they contribute a total of
$t\delta (X_{w,1,id}\otimes x)$ to the sum. Clearly case (b) makes no
contribution.

From case (c) we obtain a total contribution of 
\begin{equation}\label{ceq}
\sum_{\substack{1\leq i\leq r\\ i\ \text{free in } w\ }}
\sum_{\substack{1\leq m\leq r\\ m\ \text{joined in } w}}
(i,m)(X_{w,1,id}\otimes x)\  +
\sum_{\substack{r< j\leq r+s\\ j\ \text{free in } w\ }}
\sum_{\substack{r< l\leq r+s\\ l\ \text{joined in } w}}
 (j,l)(X_{w,1,id}\otimes x).
\end{equation}

From case (d), note that each contribution occurs twice, as it comes
from the action of $E_{i,j}$ and of $E_{l,m}$. Therefore the total
contribution can be obtained by summing over all pairs which are not
free to the left of the wall, together with the corresponding pairs to
the right of the wall. This corresponds to 
\begin{multline}\label{deq}
\Bigl(\sum_{1\leq i<m\leq r}(i,m)\ +\sum_{r<j<l\leq r+s}(j,l)
\ -\sum_{\substack{1\leq i<m\leq r\\ i,m\ \text{free in } w\ }}
(i,m)\ -
\sum_{\substack{r< j<l\leq r+s\\ j,l\ \text{free in } w\ }}
 (j,l)\\
\ -
\sum_{\substack{1\leq i\leq r\\ i\ \text{free in } w\ }}
\sum_{\substack{1\leq m\leq r\\ m\ \text{joined in } w}}
(i,m)\  -
\sum_{\substack{r< j\leq r+s\\ j\ \text{free in } w\ }}
\sum_{\substack{r< l\leq r+s\\ l\ \text{joined in } w}}
 (j,l)
\Bigr)(X_{w,1,id}\otimes x).$$
\end{multline}

Adding the four cases we see that the two double summations in
\eqref{deq} are cancelled out by the terms in \eqref{ceq}, and hence
we see that 
\begin{multline*}
T_{r,s}(X_{w,1,id}\otimes x)=\\
\Bigl(t\delta+
\sum_{1\leq i<m\leq r}(i,m)\ +\sum_{r<j<l\leq r+s}(j,l)
\ -\sum_{\substack{1\leq i<m\leq r\\ i,m\ \text{free in } w\ }}
(i,m)\ -
\sum_{\substack{r< j<l\leq r+s\\ j,l\ \text{free in } w\ }}
 (j,l)\Bigr)(X_{w,1,id}\otimes x).
\end{multline*}

By our identifications, the third sum in this equation corresponds to
the standard action of $\sum_{1\leq i<m\leq r}(i,m)$ in $\Sigma_r$ on
the first Specht module $S^{\lambda^L}$. But this sum is a central
element in $\Sigma_r$ acting on an irreducible module, and hence acts
as a scalar. By \cite[Chapter 1]{diaconis} this scalar is given by
$\sum_{d\in[\lambda^L]}c(d)$. Similarly the final sum correspond to the
action of the sum of all transpositions in $\Sigma_s$ on $S^{\lambda^R}$,
and hence acts as the scalar $\sum_{d\in[\lambda^R]}c(d)$. Substituting for
these scalars gives the desired result when $k=\CC$.

For general $k$ note that the cell
modules are all defined over $\ZZ[\delta]$. Thus our result must be
true over $\ZZ[\delta]$, and hence by base change over any field.
\end{proof}

\begin{thm} \label{deltacond}
Suppose that
$[\Delta_{r,s}(\mu^L,\mu^R):L_{r,s}(\lambda^L,\lambda^R)]\neq 0$.
Then either $(\lambda^L,\lambda^R)= (\mu^L,\mu^R)$ or
$(\lambda^L,\lambda^R)\in\Lambda^{r-a,s-a}$ and
$(\mu^L,\mu^R)\in\Lambda^{r-b,s-b}$ for some $b-a=t\geq 0$. Further we must
have that
$$t\delta+\sum_{d\in[\lambda^L]}c(d)+\sum_{d\in[\lambda^R]}c(d)
-\sum_{d\in[\mu^L]}c(d)-\sum_{d\in[\mu^R]}c(d)=0.
$$
\end{thm}
\begin{proof}
A slightly stronger version of the first part of the theorem (with
$t>0$) has already been noted in the quasi-hereditary case. The weaker
version here holds in all cases by the cellular structure of
$B_{r,s}$. For the second part, it follows from the exactness of the
localisation functor that
$$[\Delta_{r,s}(\mu^L,\mu^R):L_{r,s}(\lambda^L,\lambda^R)]=
[\Delta_{r-a,s-a}(\mu^L,\mu^R):L_{r-a,s-a}(\lambda^L,\lambda^R)]$$ and
hence we may assume that $(\lambda^L,\lambda^R)\in\Lambda^{r,s}$. If
$k=\CC$ we have that
$L_{r,s}(\lambda^L,\lambda^R)=\Delta_{r,s}(\lambda^L,\lambda^R) \cong
S^{\lambda^L}\boxtimes S^{\lambda^R}$, the lift of the irreducible for
$k\Sigma_{r,s}$, and so any walled Brauer diagram having fewer than
$r+s$ propagating lines must act as zero. In particular this includes
the action of the element $T_{r,s}$. By base change via $\ZZ[\delta]$
this must hold over any field.

By assumption there exists a $B_{r,s}$-submodule $M\leq
\Delta_{r,s}(\mu^L,\mu^R)$ and a $B_{r,s}$-homomorphism 
$$\phi: L_{r,s}(\lambda^L,\lambda^R)\longrightarrow
\Delta_{r,s}(\mu^L,\mu^R)/M.$$ 

If $k=\CC$ then the action of $$\sum_{1\leq i<m\leq r}(i,m)+ \sum_{r<
i<m\leq r+s}(i,m)$$ in the centre of $\CC\Sigma_{r,s}$ on
$S^{\lambda^L}\boxtimes S^{\lambda^R}$ must be by a scalar, and by another
application of \cite[Chapter 1]{diaconis} this equals
$$\sum_{d\in[\lambda^L]}c(d)+\sum_{d\in[\lambda^R]}c(d).$$ The same is
true for general fields exactly as before. Hence
$$\Bigl(\sum_{1\leq i<m\leq r}(i,m)+ \sum_{r< i<m\leq
r+s}(i,m)\Bigr)\phi(x)=
\Bigl(\sum_{d\in[\lambda^L]}c(d)+\sum_{d\in[\lambda^R]}c(d)\Bigr)\phi(x)$$
for all $x\in\Delta_{r,s}(\lambda^L,\lambda^R)$, and so for all
$y+M\in\im\phi$ we must have
$$T_{r,s}(y+M)=\Bigl(t\delta+\sum_{d\in[\lambda^L]}c(d)
+\sum_{d\in[\lambda^R]}c(d)-
\sum_{d\in[\mu^L]}c(d)-\sum_{d\in[\mu^R]}c(d)\Bigr)(y+M)$$ by Lemma
\ref{keyscalar}. But $T_{r,s}$ must act by zero, and so the result
follows.
\end{proof}

By standard cellular arguments
\cite[(3.9.8)]{gl} we deduce

\begin{cor}\label{weakdelta}
 Suppose that $(\lambda^L,\lambda^R)\in\Lambda^{r-a,s-a}$ and
$(\mu^L,\mu^R)\in\Lambda^{r-b,s-b}$ for some $b-a=t\geq 0$. If
$L_{r,s}(\lambda^L,\lambda^R)$ and $L_{r,s}(\mu^L,\mu^R)$ are in the same
block then
$$t\delta+\sum_{d\in[\lambda^L]}c(d)+\sum_{d\in[\lambda^R]}c(d)
-\sum_{d\in[\mu^L]}c(d)-\sum_{d\in[\mu^R]}c(d)=0.
$$
\end{cor}

We are now able to give a complete description of the blocks of
$B_{r,s}$ when $\delta$ in not an integer.

\begin{thm}\label{notintcase}
Suppose that $\delta\notin \ZZ$ and $k$ is arbitrary. Then two simple
$B_{r,s}$ modules $L(\lambda^L,\lambda^R)$ and $L(\mu^L,\mu^R)$ are in
the same block if and only if $|\lambda^L|=|\mu^L|$ (and hence
$|\lambda^R|=|\mu^R|$) and the corresponding simple
$k\Sigma_{|\lambda^L|,|\lambda^R|}$-modules are in the same block.

In particular, $B_{r,s}$ is semisimple if $\delta$ is not an integer
and $k$ is $\Sigma$-semisimple.
\end{thm}
\begin{proof}
The result follows immediately from Corollary \ref{weakdelta}, and the
fact that via localisation any two modules in the same block must both
arise as lifts from the given group algebra.
\end{proof}

\section{A sufficient condition for semisimplicity}\label{sufsec}

In the preceding section we saw that if $k$ is $\Sigma$-semisimple
then $B_{r,s}$ is semisimple for $\delta\notin \ZZ$, and it also
followed that over such fields $B_{r,s}$ was semisimple for
$|\delta|>>0$. We would like to give a stronger semisimplicity
criterion, which we will achieve by refining our condition for the
existence of a (nonzero) homomorphism between cell modules. In
this section we assume that $k$ is $\Sigma$-semisimple; clearly
$B_{r,s}$ must be non-semisimple in the other cases, by the
non-semisimplicity of $\Sigma_{r,s}$. Throughout this section we
leave to the reader the easy modifications to the proofs required
for the case $r=s$ with $\delta=0$.

Given two partitions $\lambda$ and $\mu$ we write
$\lambda\subseteq\mu$ if $\lambda$ is a subpartition of $\mu$ (i.e.~
$\lambda_i\leq \mu_i$ for all $i$).

\begin{prop}\label{mustinclude}
Suppose that $(\lambda^L,\lambda^R)\vdash (a,b)$ and
$(\mu^L,\mu^R)\vdash (a-t,b-t)$ for some $t\geq 0$ are such that
\begin{equation}\label{ifhom}
\Hom(\Delta_{r,s}(\lambda^L,\lambda^R),\Delta_{r,s}(\mu^L,\mu^R))\neq 0.
\end{equation}
Then we must have $\mu^L\subseteq \lambda^L$ and $\mu^R\subseteq
\lambda^R$ (which we will write as $(\mu^L,\mu^R)\subseteq
(\lambda^L,\lambda^R)$).
\end{prop}

\begin{proof}
 By Theorem
\ref{reduceto}(i) we may assume that $(a,b)=(r,s)$. We proceed by
induction on $t$ and on $r+s$. The result is clear when $t=0$ by the
quasihereditary structure of $B_{r,s}$; note that the base cases where
$r$ or $s$ equal $1$ for the induction on $r+s$ will be considered as
part of the inductive argument.

Suppose $t>0$ (and hence $r>0$ and $s>0$). By  Corollary \ref{indrule}
and $\Sigma$-semisimplicity
we have for any removable box $\square\in\lambda^L$ that
$$\Delta_{r,s}(\lambda^L,\lambda^R)\subseteq 
\hd(\ind_{r-1,s}^L\Delta_{r-1,s}(\lambda^L-\square,\lambda^R)).$$ 
Therefore our assumption (\ref{ifhom}) implies that
$$\Hom(\ind_{r-1,s}^L\Delta_{r-1,s}(\lambda^L-\square,\lambda^R),
\Delta_{r,s}(\mu^L,\mu^R))\neq 0 $$
and by Frobenius reciprocity we deduce that
\begin{equation}\label{frsays}
\Hom(\Delta_{r-1,s}(\lambda^L-\square,\lambda^R),
\res_{r,s}^L\Delta_{r,s}(\mu^L,\mu^R))\neq 0.
\end{equation}
By Theorem \ref{resrule} we have a short exact sequence
$$0\too\!\!\!\!
\biguplus_{\square'\in\remo(\mu^L)}\!\!\!
\Delta_{r-1,s}(\mu^L-\square',\mu^R)
\too \res_{r,s}^L\Delta_{r,s}(\mu^L,\mu^R) \too\!\!\!\!
\biguplus_{\square''\in\add(\mu^R)}\!\!\!
\Delta_{r-1,s}(\mu^L,\mu^R+\square'')
\too 0$$ 
where the left-hand term is zero if $\mu^L=\emptyset$,
and hence from (\ref{frsays}) we must have either 
\begin{equation}\label{resone}
\Hom(\Delta_{r-1,s}(\lambda^L-\square,\lambda^R),
\Delta_{r-1,s}(\mu^L-\square',\mu^R))
\neq 0
\end{equation}
for some removable box $\square'$ for $\mu^L$ or
\begin{equation}\label{restwo}
\Hom(\Delta_{r-1,s}(\lambda^L-\square,\lambda^R),
\Delta_{r-1,s}(\mu^L,\mu^R+\square''))
\neq 0
\end{equation}
for some addable box $\square''$ for $\mu^R$.

In the case (\ref{restwo}) we have that
$(\lambda^L-\square,\lambda^R)\vdash (a-1,b)$ and
$(\mu^L,\mu^R+\square'')\vdash(a-1-(t-1),b-(t-1))$ and so
$\mu^L\subseteq\lambda^L-\square$ and
$\mu^R+\square''\subseteq\lambda^R$ by the inductive hypothesis on
$t$. But clearly this implies that $(\mu^L,\mu^R)\subseteq
(\lambda^L,\lambda^R)$.

In the case (\ref{resone}) (which cannot occur when $r=1$) we have
that $(\lambda^L-\square,\lambda^R)\vdash (a-1,b)$ and
$(\mu^L-\square',\mu^R)\vdash(a-1-t,b-t)$ and so
$\mu^L-\square'\subseteq\lambda^L-\square$ and
$\mu^R\subseteq\lambda^R$ by the inductive hypothesis on $r+s$. This
only implies that $\mu^R\subseteq\lambda^R$, however repeating the above
argument with $\ind^R$ and $\res^R$ instead of $\ind^L$ and $\res^L$
also gives that $\mu^L\subseteq\lambda^L$, and so we are done.
\end{proof}

\begin{cor}\label{semisimpif}
The algebra $B_{r,s}(\delta)$ is semisimple if $k$ is
$\Sigma$-semisimple and $|\delta|\geq r+s-1$.
\end{cor}

\begin{proof}
By Theorem \ref{reduceto}(ii) it is enough to show that there are no
homomorphisms from $\Delta_{r,s}(\lambda^L,\lambda^R)$ to
$\Delta_{r,s}(\mu^L,\mu^R)$ with $(\lambda^L,\lambda^R)\vdash(r,s)$
and $(\mu^L,\mu^R)\vdash (r-1,s-1)$. By Proposition \ref{mustinclude},
the existence of such a homomorphism implies that
$\lambda^L=\mu^L+\square$ and $\lambda^R=\mu^R+\square'$ for some
addable boxes $\square$ and $\square'$.

As $\lambda^L\vdash r$ and $\lambda^R\vdash s$ we must have
$|c(\square)|\leq r-1$ and $|c(\square')|\leq s-1$. But Theorem
\ref{deltacond} implies that
$$\delta+c(\square)+c(\square')=0$$
which then gives the desired result.
\end{proof}

\section{A semisimplicity criterion}\label{twosec}

In this section we will complete the classification of semisimple
walled Brauer algebras by constructing non-zero homomorphisms between
cell modules whose weights differ by two boxes (and satisfy the
content condition in Theorem \ref{deltacond}), when $k$ is
$\Sigma$-semisimple.  These form a very special class of
homomorphisms, where the two weights are as near as is possible for
such a non-trivial homomorphism to exist, but it will turn out that
these will also be sufficient to determine the blocks in
characteristic zero.

\emph{Throughout this section we will assume (unless otherwise stated) that
$k$ is $\Sigma$-semisimple.} 

Let $c^{\lambda}_{\mu\nu}$ be the Littlewood-Richardson coefficient
denoting the multiplicity of $S^{\lambda}$ in the module
$\ind_{\CC\Sigma_{r,s}}^{\CC\Sigma_{r+s}}(S^{\mu}\boxtimes
S^{\nu})$. We will need the following result of Halverson
\cite[Corollary 7.24]{halvwall} (an analogue of \cite[Theorem
4.1]{hw3}), describing the decomposition of cell modules for
$B_{r,s}$ when regarded as $\Sigma_{r,s}$-modules, in terms of the
$c^{\lambda}_{\mu\nu}$. Although stated only for $k=\CC$ in
\cite{halvwall}, it clearly holds whenever $k$ is $\Sigma$-semisimple.

\begin{thm}\label{halvthm}
Suppose that $k$ is $\Sigma$-semisimple and
$(\mu^L,\mu^R)\in\Lambda^{r-t,s-t}$. Then
$$[\res_{k\Sigma_{r,s}}\Delta_{r,s}(\mu^L,\mu^R):S^{\lambda^L}\boxtimes
S^{\lambda^R}]=\sum_{\tau\,\vdash\, t}
c^{\lambda^L}_{\mu^L\tau}c^{\lambda^R}_{\mu^R\tau}.$$
\end{thm}

We note that Halverson's result can also be used to give an
alternative proof of Proposition \ref{mustinclude}.

We say that $(\mu^L,\mu^R)\lhd(\lambda^L,\lambda^R)$ if
$(\mu^L,\mu^R)$ can be obtained from $(\lambda^L,\lambda^R)$ by
removing a box from $\lambda^L$ and from $\lambda^R$. Then we have in
the case $t=1$ that
\begin{equation}\label{oneres}
\res_{k\Sigma_{r,s}}\Delta_{r,s}(\mu^L,\mu^R)\cong
\bigoplus_{(\mu^L,\mu^R)\lhd(\lambda^L,\lambda^R)}S^{\lambda^L}\boxtimes
S^{\lambda^R}.
\end{equation}

We will now deduce the existence of the desired two box
homomorphisms. The proof is similar to that of \cite[Theorem
5.2]{cdm}. For two partitions $\lambda$ and $\mu$, we denote by
$\lambda/\mu$ the skew partition whose skew Young diagram
$[\lambda/\mu]$ consists of those boxes in $[\lambda]$ which are not in
$[\mu]$.  Given $\mu\subset\lambda$ with $|\lambda|-|\mu|=1$ we denote by
$c(\lambda/\mu)$ the content of the unique box in 
$[\lambda/\mu]$.

\begin{thm}\label{twoboxhom}
Let $k$ be $\Sigma$-semisimple. Suppose that
$(\lambda^L,\lambda^R)\in\Lambda^{r-t,s-t}$, and
$(\mu^L,\mu^R)\in\Lambda^{r-t-1,s-t-1}$. If we have that
$(\mu^L,\mu^R)\lhd (\lambda^L,\lambda^R)$ and
$$c(\lambda^L/\mu^L)+c(\lambda^R/\mu^R)+\delta=0$$
then
$$\Hom_{B_{r,s}}(\Delta_{r,s}(\lambda^L,\lambda^R),
\Delta_{r,s}(\mu^L,\mu^R))\cong k.$$
\end{thm}
\begin{proof}
As the localisation functor is exact we may assume that
$(\lambda^L,\lambda^R)\vdash (r,s)$, and hence that
$$\Delta_{r,s}(\lambda^L,\lambda^R)\cong S^{\lambda^L}\boxtimes
S^{\lambda^R}$$
(where all Brauer diagrams with propagating vector other than $(r,s)$
act as zero on the righthand side).

We will fix a labelling of the boxes of $\lambda^L$ and
$\lambda^R$. Number the boxes of $\lambda^L$ with $1,2,\ldots,r$ along
the rows from left to right, and from top to bottom. Similarly number
the boxes of $\lambda^R$ with $r+1,\ldots,r+s$. Using the natural
identifications of $\Sigma_r$ with Sym$\{1,2,\ldots r\}$ and 
$\Sigma_s$ with Sym$\{r+1,r+2,\ldots r+s\}$ we can define elements
$e_{\lambda^L}$ and $e_{\lambda^R}$ by setting
$$e_{\lambda^L}=\frac{f^{\lambda^L}}{r!}\sum_{\sigma_L\in
  C(\lambda^L)}\sum_{\tau_L\in
  R(\lambda^L)}\sgn(\sigma_L)\sigma_L\tau_L$$
and 
$$e_{\lambda^R}=\frac{f^{\lambda^R}}{s!}\sum_{\sigma_R\in
  C(\lambda^R)}\sum_{\tau_R\in
  R(\lambda^R)}\sgn(\sigma_R)\sigma_R\tau_R$$ where $f^{\lambda}=\dim
  S^{\lambda}$ and $C(\lambda)$ and $R(\lambda)$ are respectively the
  column and row stabilisers of $\lambda$. These are idempotents, such
  that $e_{\lambda^L}\CC\Sigma_r\cong a_{\lambda}S^{\lambda^L}$ and
  $e_{\lambda^R}\CC\Sigma_s\cong b_{\lambda}S^{\lambda^R}$ where the
  multiplicities $a_\lambda$ and $b_{\lambda}$ will not concern us
  (see for example \cite[Chapter 7]{fultab}). (It is clear that the
  denominators in each idempotent are non-zero when $k$ is
  $\Sigma$-semisimple. The numerator is also non-zero as dimensions of
  Specht modules can be written as products of hook lengths, and hence
  cannot involve factors bigger than the total number of boxes.)  Note
  that the idempotent
$$e_{\lambda}=e_{\lambda^L}e_{\lambda^R}=
e_{\lambda^R}e_{\lambda^L}\in k\Sigma_{r,s}$$
is such that for any $\Sigma_{r,s}$-module $M$, $e_{\lambda}M$ is the
$S^{\lambda^L}\boxtimes S^{\lambda^R}$ isotypic component of $M$. 

Let $W=e_{\lambda}\Delta_{r,s}(\mu^L,\mu^R)$. By  (\ref{oneres}) we
know that 
$$W\cong S^{\lambda^L}\boxtimes S^{\lambda^R}$$ as a
$\Sigma_{r,s}$-module. To show that this is in fact a
$B_{r,s}$-submodule of $\Delta_{r,s}(\mu^L,\mu^R)$, it is enough to
show that $E_{i,j}W=0$ for all $1\leq i\leq r<j\leq r+s$. Indeed, it is
enough to show that this holds for a single choice of $i$ and $j$, as 
$$\sigma E_{i,j}\sigma^{-1}=E_{\sigma(i),\sigma(j)}$$
for all $\sigma\in\Sigma_{r,s}$.

Fix  $1\leq i\leq r<j\leq r+s$ and consider the map 
$$E_{i,j}:\Delta_{r,s}(\mu^L,\mu^R)\too \Delta_{r,s}(\mu^L,\mu^R).$$
This is a $k\Sigma_{r-1,s-1}$-homomorphism, where we identify 
$$\Sigma_{r-1,s-1}\cong \mbox{\rm Sym}(\{1,2,\ldots
r\}\backslash\{i\})
\times \mbox{\rm Sym}(\{r+1,r+2,\ldots r+s\}\backslash\{j\}).$$
Note that $E_{i,j}(\Delta_{r,s}(\mu^L,\mu^R))\subseteq U$ where $U$ is
the span of all elements of the form $X_{w_0,1,id}\otimes x$ where $w_0$
has an arc between $i$ and $j$ and $x\in S^{\mu^L}\boxtimes S^{\mu^R}$,
and 
\begin{equation}\label{Ures}
\res_{k\Sigma{r-1,s-1}}U\cong S^{\mu^L}\boxtimes S^{\mu^R}
\end{equation}
By (\ref{oneres}) we have that 
$$\res_{k\Sigma{r-1,s-1}}W\cong
\bigoplus_{(\nu^L,\nu^R)\lhd(\lambda^L,\lambda^R)} S^{\nu^L}\boxtimes
S^{\nu^R}$$ and by (\ref{Ures}) every summand in this direct sum must
be sent by $E_{i,j}$ to zero except possibly for $V=S^{\mu^L}\boxtimes
S^{\mu^R}$. Write
$\res_{k\Sigma{r-1,s-1}}W=V\oplus Y$.

Up until this point $i$ and $j$ have been arbitrary. Henceforth we
will consider the case where $i$ is the number labelling the box in 
$\lambda^L/\mu^L$ and $j$ is the number labelling the box in 
$\lambda^R/\mu^R$. For example, if $\lambda^L=(3,2,1)$ and
$\lambda^R=(2,2)$, with $\mu^L=(2,2,1)$ and $\mu^R=(2,1)$, then $i=3$
and $j=10$ as illustrated in Figure \ref{ijex}.

\begin{figure}[ht]
\includegraphics{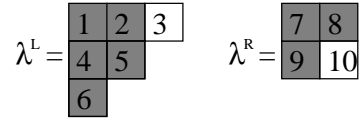}
\caption{$\lambda^L=(3,2,1)$ and $\lambda^R=(2,2)$,
  with $\mu^L=(2,2,1)$ and $\mu^R=(2,1)$ shaded.}
\label{ijex}
\end{figure}

Writing $e_{\lambda}(X_{w_0,1,id}\otimes x)=v+y$ with $v\in V$ and $y\in
Y$ (which is independent of $\delta$) the remarks above imply that
$$E_{i,j}e_{\lambda}(X_{w_0,1,id}\otimes x)=E_{i,j}v.$$

We claim that the coefficient of $X_{w_0,1,id}\otimes x$ in $E_{i,j}v$ is a
non-zero multiple of 
$$\delta+c(\lambda^L/\mu^L)+c(\lambda^R/\mu^R).$$
This would imply that $v\neq 0$, and that for 
$$\delta+c(\lambda^L/\mu^L)+c(\lambda^R/\mu^R)=0$$
we would have $E_{i,j}V=0$ (as $V$ is a simple module), and hence that
$E_{i,j}W=0$ as required. Thus it is enough to prove the claim.

We have that
$$E_{i,j}e_{\lambda}(X_{w_0,1,id}\otimes x)=
\frac{f^{\lambda^L}}{r!}\frac{f^{\lambda^R}}{s!}\!\!\!\!
\sum_{\genfrac{}{}{0cm}{1}
{\sigma_L\in C(\lambda^L)}{\sigma_R\in
C(\lambda^R)}}
\sum_{\genfrac{}{}{0cm}{1}{\tau_L\in
R(\lambda^L)}
{\tau_R\in R(\lambda^R)}}
\sgn(\sigma_L\sigma_R)E_{i,j}\sigma_L\sigma_R\tau_L\tau_R
(X_{w_0,1,id}\otimes x)$$ and we wish to find the coefficient of
$(X_{w_0,1,id}\otimes x)$ in this sum.

In order to keep track of the various cases, it will be convenient to
have a graphical notation for the partial one row diagrams arising as
configurations of northern arcs in the summands of this expression. We
will represent elements of ${\mathcal V}_{r,s,t}$ by adding ties to
the double Young tableau joining each pair of nodes connected by an
arc, omitting any such ties which do not play a role in the
calculation. For example, the element $w_0$ can be represented by the
diagram in Figure \ref{w0is}.

\begin{figure}[ht]
\includegraphics{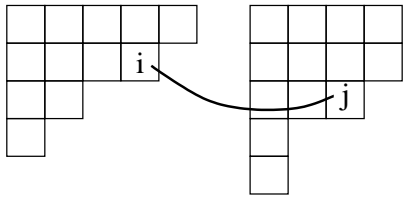}
\caption{A graphical representation of the arc $ij$ in $w_0$}
\label{w0is}
\end{figure}

\noindent
{\bf Case 1:} Suppose that $\sigma_L\sigma_R\tau_L\tau_RX_{w_0,1,id}$
has an edge between $i$ and $j$. This occurs if and only if all of
$\sigma_L$, $\sigma_R$, $\tau_L$ and $\tau_R$ fix $i$ and $j$.
 In this case
$$E_{i,j}\sigma_L\sigma_R\tau_L\tau_R(X_{w_0,1,id}\otimes x)
=\delta\sigma_L\sigma_R\tau_L\tau_R(X_{w_0,1,id}\otimes x).$$ 
For $\sigma_L\sigma_R\tau_L\tau_R(X_{w_0,1,id}\otimes x)$ to be in the
span of $X_{w_0,1,id}\otimes S^{\lambda^L}\boxtimes S^{\lambda^R}$ we
must have
$$\begin{array}{l} \tau_L\in R(\mu^L)\subset R(\lambda^L)\\ 
\sigma_L\in C(\mu^L)\subset C(\lambda^L)\\ 
\tau_R\in R(\mu^R)\subset R(\lambda^R)\\
\sigma_R\in C(\mu^R)\subset C(\lambda^R).\end{array}$$ 
As any such quartet of elements fixes $i$ and
$j$, the contribution to our sum in this case equals
$$\frac{f^{\lambda^L}}{r!}\frac{f^{\lambda^R}}{s!}\!\!\!\!
\sum_{\genfrac{}{}{0cm}{1} {\sigma_L\in
C(\mu^L)}{\sigma_R\in C(\mu^R)}}
\sum_{\genfrac{}{}{0cm}{1}{\tau_L\in R(\mu^L)}
{\tau_R\in R(\mu^R)}}
\delta\sgn(\sigma_L\sigma_R)E_{i,j}\sigma_L\sigma_R\tau_L\tau_R
(X_{w_0,1,id}\otimes x)$$ which  equals
$$\delta\frac{f^{\lambda^L}}{r!}\frac{f^{\lambda^R}}{s!}
\frac{(r-1)!}{f^{\mu^L}}\frac{(s-1)!}{f^{\mu^R}}
X_{w_0,1,id}\otimes e_{\mu}x.$$
Now $e_{\mu}(x)=x$ for all $x\in
S^{\mu^L}\boxtimes S^{\mu^R}$, and hence we obtain the contribution 
$$\delta\frac{f^{\lambda^L}}{f^{\mu^L}}\frac{f^{\lambda^R}}{f^{\mu^R}}
\frac{1}{rs}
X_{w_0,1,id}\otimes x$$
from the terms arising in Case 1.

\noindent
{\bf Case 2:} Suppose that neither $i$ nor $j$ is part of an arc
in $\sigma_L\sigma_R\tau_L\tau_RX_{w_0,1,id}$. Then
$$E_{i,j}\sigma_L\sigma_R\tau_L\tau_R(X_{w_0,1,id}\otimes x)=0$$ 
so such terms contribute nothing to our sum.

\noindent
{\bf Case 3:} The only remaining case is when exactly one of $i$ or
$j$ is part of an arc in
$\sigma_L\sigma_R\tau_L\tau_RX_{w_0,1,id}$. (Note that they cannot
both be, as there is only one arc in such a term.) We will consider
the subcase where the arc connects $i$ to some element $k$ with $r+1\leq
k\leq r+s$ (the subcase of an edge between $j$ and some element $k$ with
$1\leq k\leq r$ is similar). We must have
\begin{equation}\label{3aorb1}
\sigma_L\in C(\mu^L)\quad\quad \wand\quad\quad
\tau_L\in R(\mu^L)
\end{equation} 
There are two possibilities: (a) $k$ is in the same column as $j$, or (b)
$k$ is in a column to the left of the column containing $j$

\noindent
{\bf Subcase 3(a):} Suppose that $k$ is in the same column as $j$, as
illustrated in Figure \ref{case3a}.

\begin{figure}[ht]
\includegraphics{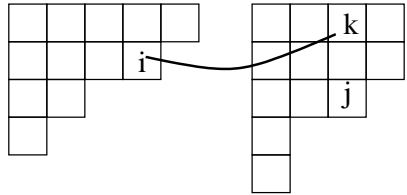}
\caption{The subcase 3(a)}
\label{case3a}
\end{figure}

We must have 
\begin{equation}\label{3a2}
\tau_R\in R(\mu^R)\quad\quad\wand
\quad\quad \sigma_R\in(k,j)C(\mu^R)
\end{equation}
and every quartet satisfying (\ref{3aorb1}) and (\ref{3a2}) will arise
in this way. Writing $\sigma_R=(k,j)\sigma_R'$ with $\sigma_R'\in
C(\mu^R)$, and noting that
$\sgn(\sigma_R')=-\sgn(\sigma_R)$, we see that we obtain a
contribution of
$$\sum_{k \ \mbox{\rm\small  above}\ j}
\frac{f^{\lambda^L}}{r!}\frac{f^{\lambda^R}}{s!}\!\!\!\!
\sum_{\genfrac{}{}{0cm}{1} {\sigma_L\in
C(\mu^L)}{\sigma_R'\in C(\mu^R)}}
\sum_{\genfrac{}{}{0cm}{1}{\tau_L\in R(\mu^L)}
{\tau_R\in R(\mu^R)}}
(-\sgn(\sigma_L\sigma_R'))E_{i,j}\sigma_L(k,j)\sigma_R'\tau_L\tau_R
(X_{w_0,1,id}\otimes x).$$
Arguing as in Case 1 we see that this equals
$$-\sum_{k \ \mbox{\rm\small above}\ j}
\frac{f^{\lambda^L}}{f^{\mu^L}} \frac{f^{\lambda^R}}{f^{\mu^R}}
\frac{1}{rs}E_{i,j}(k,j) X_{w_0,1,id}\otimes e_{\mu}x
=-\abv_{\lambda^R}(j) \frac{f^{\lambda^L}}{f^{\mu^L}}
\frac{f^{\lambda^R}}{f^{\mu^R}} \frac{1}{rs}X_{w_0,1,id}\otimes x$$
where $\abv_{\lambda^R}(j)$ denotes the number of boxes above $j$ in
$\lambda^R$.

\noindent
{\bf Subcase 3(b):} Suppose that $k$ is in a column to the left of the
column containing $j$, with $l$ as illustrated in Figure \ref{case3b}.

\begin{figure}[ht]
\includegraphics{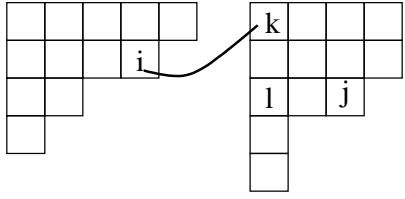}
\caption{The subcase 3(b)}
\label{case3b}
\end{figure}

We must have 
\begin{equation}\label{3b2}
\tau_R\in (j,l)R(\mu^R)\quad\quad\wand
\quad\quad \sigma_R\in C(\mu^R)
\end{equation}
(as $k$ is an arbitrary element in the column containing $l$) and
every quartet satisfying (\ref{3aorb1}) and (\ref{3b2}) will arise in
this way. Arguing as in Subcase 3(a), but now with no sgn
modifications, we see that we obtain a contribution of 
$$\lft_{\lambda^R}(j) \frac{f^{\lambda^L}}{f^{\mu^L}}
\frac{f^{\lambda^R}}{f^{\mu^R}} \frac{1}{rs}X_{w_0,1,id}\otimes x$$
where $\lft_{\lambda^R}(j)$ denotes the number of boxes to the left of 
$j$ in $\lambda^R$.

Combining Subcase 3(a) and 3(b), and the corresponding versions where
the arc connects to $j$ rather than $i$, we obtain a total
contribution in Case 3 of
$$\big(\lft_{\lambda^R}(j) -\abv_{\lambda^R}(j) 
+\lft_{\lambda^L}(i) -\abv_{\lambda^L}(i)\big) 
\frac{f^{\lambda^L}}{f^{\mu^L}}
\frac{f^{\lambda^R}}{f^{\mu^R}} \frac{1}{rs}
X_{w_0,1,id}\otimes x.$$
As
$$\lft_{\lambda^R}(j) -\abv_{\lambda^R}(j) +\lft_{\lambda^L}(i)
-\abv_{\lambda^L}(i)=c(\lambda^L/\mu^L)+c(\lambda^R/\mu^R)$$
combining Cases 1--3 now gives that the coefficient of 
$X_{w_0,1,id}\otimes x$ in $E_{i,j}v$ equals
$$\frac{f^{\lambda^L}}{f^{\mu^L}}
\frac{f^{\lambda^R}}{f^{\mu^R}} \frac{1}{rs}(\delta+
c(\lambda^L/\mu^L)+c(\lambda^R/\mu^R))$$
which is a nonzero multiple of 
$\delta+ c(\lambda^L/\mu^L)+c(\lambda^R/\mu^R)$
as required.
\end{proof}

Combining the last result with the semisimplicity results in Section
\ref{sufsec}  we obtain

\begin{thm}\label{whenss}
The walled Brauer algebra $B_{r,s}(\delta)$ is semisimple
if and only if $k$ is $\Sigma$-semisimple and 
one of the following conditions holds:\\
(i) $\ \delta\notin\ZZ$, or\\ 
(ii) $|\delta|> r+s-2$, or\\ 
(iii) $r=0$ or $s=0$, or\\
(iii) $\delta=0$ and $(r,s)=(1,2)$, $(1,3)$, $(2,1)$, or $(3,1)$.
\end{thm}
\begin{proof}
Clearly $B_{r,s}$ cannot be semisimple if $k$ is not
$\Sigma$-semisimple. If $k$ is $\Sigma$-semisimple and $r$ or $s$
equals zero then $B_{r,s}$ must be semisimple (by the definition of
$\Sigma$-semisimplicity as it is just the group algebra of the
symmetric group $\Sigma_{r,s}$). By Theorem \ref{notintcase} and Corollary
\ref{semisimpif} the only other cases where $B_{r,s}$ can be
non-semisimple occur when $\delta\in\ZZ$ and $|\delta|\leq r+s-2$.  By
Theorem \ref{reduceto} it is enough to classify exactly those pairs
$(r,s)$ for which there exists $(\lambda^L,\lambda^R)\vdash (r,s)$ and
$(\mu^L,\mu^R)\vdash (r-1,s-1)$ with a non-zero homomorphism from
$\Delta_{r,s}(\lambda^L,\lambda^R)$ to $\Delta_{r,s}(\mu^L,\mu^R)$.

We may assume that $\delta\in\ZZ$ with $|\delta|\leq r+s-2$. We write
$\delta=\delta_0+\delta_1$ with $|\delta_0|\leq r-1$ and
$|\delta_1|\leq s-1$. First suppose that we can choose such a
decomposition with both $\delta_0$ and $\delta_1$ non-zero.  Then we
can construct explicitly a quartet $(\lambda^L,\lambda^R)$ and
$(\mu^L,\mu^R)$ with a nonzero homomorphism in the following
manner. If $\delta_0> 0$ set
$$\lambda^L=((\delta_0+1),1^{r-\delta_0-1})\supset
\mu^L=(\delta_0,1^{r-\delta_0-1})$$ and if $\delta_0<0$ take the
transpose of this pair for $-\delta_0$. Similarly construct
$\lambda^R$ and $\mu^R$ in terms of $\delta_1$. In each case the
resulting bipartitions satisfy the conditions of Theorem
\ref{twoboxhom} and so we have a non-zero homomorphism.

The only cases where we are forced to take $\delta_0$ or $\delta_1$
equal to zero occur when $r=1$ or $s=1$. We consider the case $s=1$,
the other is similar. Clearly $\lambda^R=(1)$ has a removable box of
content zero, so we just have to determine when there exists
$\lambda^L\vdash r$ with a removable box of content $-\delta$. If
$\delta\neq 0$ this is always possible: if $\delta>0$ then we take
$\lambda^L= (r-\delta,1^{\delta-1})$, while for $\delta<0$ we take the
transpose of this partition.

Finally we are left with the case where $\delta=0$, and so require
$\lambda^L\vdash r$ with a removable box of content $0$.  Such a box
exists when $\lambda^L=(1)$ or $\lambda^L=(2,2,1^{r-4})$, but not when
$r=2$ or $r=3$. This (and the corresponding argument when $r=1$)
provides the exceptional semisimple cases listed above.
\end{proof}

\section{The blocks of the walled Brauer algebra in characteristic
  zero}\label{blocksec}

In this section we will determine the blocks of the walled Brauer in
the case when $k$ is $\Sigma$-semisimple.  Our approach is
modelled on that for the corresponding result for the Brauer algebra
in \cite{cdm}.

We begin by giving a refinement of Theorem
\ref{deltacond}. 

\begin{prop}\label{refcomp} Let $k$ be $\Sigma$-semisimple.
Suppose that $[\Delta_{r,s}(\mu^L,\mu^R):L_{r,s}(\lambda^L,\lambda^R)]\neq
0$. Then $(\mu^L,\mu^R)\subseteq (\lambda^L,\lambda^R)$, and
there exists a pairing of the boxes in $\lambda^L/\mu^L$ with those in
$\lambda^R/\mu^R$ such that the sum of the contents of the boxes in
each pair equals $-\delta$ in $k$.
\end{prop}
\begin{proof}
First note that the condition $(\mu^L,\mu^R)\subseteq
(\lambda^L,\lambda^R)$ follows immediately from Theorem \ref{halvthm}.
If $r=0$ or $s=0$ then $B_{r,s}$ is just the group algebra of the
symmetric group $\sigma_{r+s}$, and the result follows from the
definition of $\Sigma$-semisimplicity. 

We will proceed by induction on $r+s$. The case $r=s=1$ follows from
Corollary \ref{weakdelta}, as the only allowable bipartitions are
$(1,1)$ and $(0,0)$, and hence the result is true if $r+s=2$.  Thus we
assume that the result is true for all $B_{a,b}$ with $a+b=n-1$ and
will show that it is also true for $B_{r,s}$ with $r+s=n$ and $r,s\neq
0$.

If $[\Delta_{r,s}(\mu^L,\mu^R):L_{r,s}(\lambda^L,\lambda^R)]\neq 0$
  then by the above remarks and Corollary \ref{weakdelta} we must have
  $(\mu^L,\mu^R)\subseteq (\lambda^L,\lambda^R)$ and
\begin{equation}\label{comparewith}
t\delta+\sum_{d\in[\lambda^L/\mu^L]}c(d)+\sum_{d\in[\lambda^R/\mu^R]}c(d)=0
\end{equation}
  where $t=|\lambda^L-\mu^L|=|\lambda^R-\mu^R|$. By localising we may
  assume that $(\lambda,\mu)\vdash (r,s)$, so that
  $L_{r,s}(\lambda^L,\lambda^R)=\Delta_{r,s}(\lambda^L,\lambda^R)$.
  Thus $\Delta_{r,s}(\mu^L,\mu^R)$ has a submodule $M$ such that there
  is an injection 
$$\Delta_{r,s}(\lambda^L,\lambda^R)\hookrightarrow
\Delta_{r,s}(\mu^L,\mu^R)/M.$$ By our assumption that $r$ is nonzero
there exists a removable box $\square$ in $\lambda^L$, and by
Corollary \ref{indrule} and $\Sigma$-semisimplicity there exists a surjection
$$\ind^L_{r-1,s}\Delta_{r-1,s}(\lambda^L-\square,\lambda^R)\too
\Delta_{r,s}(\lambda^L,\lambda^R).$$
Hence we have
$$\Hom(\ind^L_{r-1,s}\Delta_{r-1,s}(\lambda^L-\square,\lambda^R),
\Delta_{r,s}(\mu^L,\mu^R)/M)\neq 0$$
and so by Frobenius reciprocity we have
$$\Hom(\Delta_{r-1,s}(\lambda^L-\square,\lambda^R),
\res^L_{r,s}\left(\Delta_{r,s}(\mu^L,\mu^R)/M\right))\neq 0.$$
This implies that 
$$L_{r-1,s}(\lambda^L-\square,\lambda^R)=
\Delta_{r-1,s}(\lambda^L-\square,\lambda^R)$$
is a composition factor of $\res^L_{r,s}(\Delta_{r,s}(\mu^L,\mu^R)$.
By Theorem \ref{resrule} we see that either (i)
$$[\Delta_{r-1,s}(\mu^L-\square',\mu^R):
L_{r-1,s}(\lambda^L-\square,\lambda^R)]\neq 0$$ 
for some $\square'\in \remo(\mu^L)$, or (ii)
$$[\Delta_{r-1,s}(\mu^L,\mu^R+\square'):
L_{r-1,s}(\lambda^L-\square,\lambda^R)]\neq 0$$ 
for some $\square'\in \add(\mu^L)$. We consider each case in turn.

In case (i), our inductive hypothesis implies that
$\mu^L-\square'\subseteq \lambda^L-\square$ and 
$$t\delta+\Big(\sum_{d\in[\lambda^L/\mu^L]}c(d)
+\sum_{d\in[\lambda^R/\mu^R]}c(d)\Big)-c(\square)+c(\square')=0.$$
Comparing with (\ref{comparewith}) we see that
$c(\square)=c(\square')$.  By induction we know that there is a pairing of
the boxes in $(\lambda^L-\square)/(\mu^L-\square')$ with those in
$\lambda^R/\mu^R$ such the contents of each pair sum to $-\delta$.
But as multisets, the set of contents in
$\lambda^L-\square/\mu^L-\square'$ and in $\lambda^L/\mu^L$ are equal,
and hence there is such a pairing between $\lambda^L/\mu^L$ and
$\lambda^R/\mu^R$ as required.

Next we consider case (ii). By induction we must have that 
$\mu^L\subseteq \lambda^L-\square$ and 
$\mu^R+\square'\subseteq \lambda^R$ with 
$$(t-1)\delta+\Big(\sum_{d\in[\lambda^L/\mu^L]}c(d)
+\sum_{d\in[\lambda^R/\mu^R]}c(d)\Big)-c(\square)-c(\square')=0.$$
Comparing with (\ref{comparewith}) we see that
$c(\square)+c(\square')=-\delta$, and by induction we know that there
is a pairing of the boxes in $(\lambda^L-\square)/\mu^L$ with those in
$\lambda^R/(\mu^R+\square')$ such the contents of each pair sum to
$-\delta$.  But then extending this pairing to one between
$\lambda^L/\mu^L$ and $\lambda^R/\mu^R$ by adding the paired boxes $\square$
and $\square'$ gives the desired result.
\end{proof}

Given two partitions $\lambda$ and $\mu$, we denote by
$\lambda\cap\mu$ the partition whose corresponding Young diagram is
the intersection of those for $\lambda$ and $\mu$.

\begin{defn}
We will say that $(\lambda^L,\lambda^R)$ and $(\mu^L,\mu^R)$ are
\emph{$\delta$-balanced} (or just \emph{balanced} when this will not
cause confusion) if there is a pairing of the boxes in
$[\lambda^L/(\lambda^L\cap\mu^L)]$ with those in
$[\lambda^R/(\lambda^R\cap\mu^R)]$ and of the boxes in
$[\mu^L/(\lambda^L\cap\mu^L)]$ with those in
$[\mu^R/(\lambda^R\cap\mu^R)]$ such that the contents of each pair sum
to $-\delta$ in $k$.
\end{defn}

Just as for Corollary \ref{weakdelta}, we deduce from Proposition
\ref{refcomp} the following partial block result.

\begin{cor}\label{oneway} 
Let $k$ be $\Sigma$-semisimple. If $(\lambda^L,\lambda^R)$ and
$(\mu^L,\mu^R)$ are in the same block for $B_{r,s}(\delta)$ then they
are $\delta$-balanced.
\end{cor}

We will show that this is in fact a necessary and sufficient condition
for block membership when $k$ is $\Sigma$-semisimple. 
Given a partition $\mu\subset\lambda$,
we denote by $\remo(\lambda/\mu)$ the set of boxes in $\remo(\lambda)$
which are not in $\mu$.

\begin{defn} Suppose that $(\mu^L,\mu^R)\subseteq(\lambda^L,\lambda^R)$ is 
a balanced pair. For each $\square_i\in\remo(\lambda^L/\mu^L)$ we wish
to consider $(\mu^L,\mu^R)^i$, the \emph{$i$-maximal balanced
sub-bipartition between $(\lambda^L,\lambda^R)$ and
$(\mu^L,\mu^R)$}. This is the maximal bipartition in
$(\lambda^L,\lambda^R)$ not containing $\square_i$ such that
$(\lambda^L,\lambda^R)$ and $(\mu^L,\mu^R)^i$ is $\delta$-balanced.
\end{defn}

We can give an explicit recursive construction of
$(\mu^L,\mu^R)^i$. Given two boxes $\square$ and $\square'$ with the
same content in a partition $\lambda$ we will say that $\square$ is
\emph{larger} than $\square'$ if $\square$ appears on a later row than
$\square'$. Suppose that $(\mu^L,\mu^R)\subseteq(\lambda^L,\lambda^R)$
is a balanced pair with $\square_i\in\remo(\lambda^L/\mu^L)$. By the
balanced pair condition there exists a largest box $\square_i'$ in
$\lambda^R/\mu^R$ such that $c(\square_i)+c(\square_i')=-\delta$. Let
$[(\lambda^L,\lambda^R)/(\mu^L,\mu^R)]_0=\{\square_i,\square_i'\}$. Given 
$[(\lambda^L,\lambda^R)/(\mu^L,\mu^R)]_m$ we set 
$$[(\lambda^L,\lambda^R)/(\mu^L,\mu^R)]_{m+1}=
[(\lambda^L,\lambda^R)/(\mu^L,\mu^R)]_m\cup A_{m+1}\cup A_{m+1}'$$
where $A_{m+1}$ is the set of boxes in $\lambda^L$ or $\lambda^R$
which are to the right of or below a box in
$[(\lambda^L,\lambda^R)/(\mu^L,\mu^R)]_m$, and $A'_{m+1}$ is the set
of boxes $\square'$ in $(\lambda^L/\mu^L,\lambda^R/\mu^R)$ whose
content satisfies $c(\square)+c(\square')=-\delta$ for some
$\square\in A_{m+1}$ with $\square$ and $\square'$ not both in the
same partition which are largest with such content. Let
$(\mu^L,\mu^R)_m^i$ be the sub-bipartition of $(\lambda^L,\lambda^R)$ with 
complement $[(\lambda^L,\lambda^R)/(\mu^L,\mu^R)]_m$.

This iterative process will eventually stabilise to produce a
$\delta$-balanced sub-bipartition $(\mu^L,\mu^R)^i$ of
$(\lambda^L,\lambda^R)$, obtained by removing a strip of boxes one box
wide from the edges of each of $\lambda^L$ and $\lambda^R$. To see
this, first note that the construction of each
$[(\lambda^L,\lambda^R)/(\mu^L,\mu^R)]_m$ clearly only involves boxes
from the edges of $\lambda^L$ and $\lambda^R$, and so produces a strip
in each at most one box wide (as each box in a given strip has
different content, and is the largest with such). Second, the only way
in which the process could terminate without producing a balanced
sub-bipartition would be if one or other of the two strips ended with
one of the boxes at the end of the first row or first column of
$\lambda^L$ or $\lambda^R$, without the other strip being
removable. But this would contradict the fact that
$(\lambda^L,\lambda^R)$ and $(\mu^L,\mu^R)$ are $\delta$-balanced. It
is also easy to see that $(\mu^L,\mu^R)^i$ is maximal in
$(\lambda^L,\lambda^R)$ with this property.

\begin{example} We will illustrate the above construction with an 
example. Let $(\lambda^L,\lambda^R)=((4^3,1^3),(5^2,2^3))$ and
$(\mu^L,\mu^R)=((2,1),(4))$ as in Figure \ref{maxex}. These form a
balanced pair with $\delta=+1$.

\begin{figure}[ht]
\includegraphics{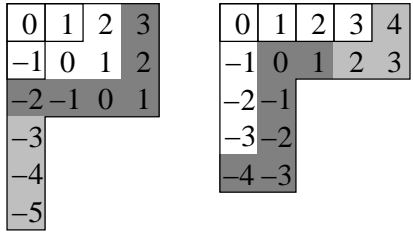}
\caption{Two examples of the $(\mu^L,\mu^R)^i$ construction.}
\label{maxex}
\end{figure}

If $\square_i$ is the largest box in $\lambda^L$ with content $-5$,
then the associated maximal balanced sub-bipartition is obtained by
removing the three lightly shaded boxes in each partition. If
$\square_i$ is the other removable box in $\lambda^L$ (with content
$1$) then the associated sub-bipartition is obtained by removing \emph{all}
shaded boxes from each partition.
\end{example}

As the above example illustrates, some of the removable strips so far
constructed may contain others. Partially order the removable strips
obtained from the above construction by inclusion. Then we define a
\emph{maximal balanced sub-bipartition $(\mu^L_{\lambda},\mu^R_{\lambda})$
between $(\lambda^L,\lambda^R)$ and $(\mu^L,\mu^R)$} to be any
balanced sub-bipartition for which the associated removable strip is
minimal. Thus in the above example there is a unique choice of
$(\mu^L_{\lambda},\mu^R_{\lambda})=((4^3),(3^2,2^3))$ given by removing the
three lightly shaded boxes. In general the choice will not be unique.

Note that if $k$ is $\Sigma$-semisimple then for each equivalence
class of integers mod $p$ there is at most one member of that class
which occurs as a content in a partition $\lambda$ of $r$ or $s$. Thus
throughout the following proof, it is unambiguous to regard all
contents mod $p$.

\begin{thm} \label{maxhom} Let $k$ be $\Sigma$-semisimple. 
If $(\mu^L,\mu^R)\subset(\lambda^L,\lambda^R)$ is a 
balanced pair then for any maximal balanced sub-bipartition
$(\mu^L_{\lambda},\mu^R_{\lambda})$ we have
$$\Hom(\Delta_{r,s}(\lambda^L,\lambda^R),
\Delta_{r,s}(\mu^L_{\lambda},\mu^R_{\lambda})))\neq 0.$$
\end{thm}

\begin{proof}
As usual, we may assume that $(\lambda^L,\lambda^R)\vdash (r,s)$. If
$(\mu^L_{\lambda},\mu^R_{\lambda})\vdash (r-1,s-1)$ then we are done by Theorem
\ref{twoboxhom}. For the remaining cases, pick $\square\in
\remo(\lambda^L/\mu^L_{\lambda})$ with $|c(\square)+\frac{\delta}{2}|$
maximal, and suppose that $\square'$ is the corresponding box in
$\lambda^R/\mu^R_{\lambda}$ with
$c(\square)+c(\square')=-\delta$.  

Note that by the maximality of $\square$ and the construction of
$(\mu^L_{\lambda},\mu^R_{\lambda})$ there is no box of content
$c(\square)$ in $\remo(\mu^L_{\lambda})$, and exactly one box of
content $c(\square')$ in $\add(\mu^R_{\lambda})$ (namely $\square'$
itself). Clearly $(\lambda^L-\square,\lambda^R)$ and
$(\mu^L_{\lambda},\mu^R_{\lambda}+\square')$ is a balanced pair; we
claim that in fact $(\mu^L_{\lambda},\mu^R_{\lambda}+\square')$ is
also a maximal balanced sub-bipartition for this balanced pair. But
this is also obvious, as any larger balanced sub-bipartition would
give rise to a corresponding balanced sub-bipartition between
$(\mu^L,\mu^R)$ and $(\lambda^L,\lambda^R)$, which would contradict
the maximality of $(\mu^L_{\lambda},\mu^R_{\lambda})$.

If $\square'\in \remo(\lambda^R)$ then
$(\mu^L_{\lambda},\mu^R_{\lambda})=
(\lambda^L-\square,\lambda^R-\square')$ by maximality, and so we are
done. Otherwise by our initial assumptions there is no
removable box in $\lambda^R$ with content $c(\square')$. Therefore by
Frobenius reciprocity, Corollary \ref{oneway}, Corollary
\ref{indrule}, and $\Sigma$-semisimplicity  we have
$$\begin{array}{ll} \Hom(\Delta_{r,s}(\lambda^L,\lambda^R),
\Delta_{r,s}(\mu^L_{\lambda},\mu^R_{\lambda}))&\cong \Hom
(\ind_{r-1,s}^L\Delta_{r-1,s}(\lambda^L-\square,\lambda^R),
\Delta_{r,s}(\mu^L_{\lambda},\mu^R_{\lambda}))\\ 
&\cong \Hom
(\Delta_{r-1,s}(\lambda^L-\square,\lambda^R), 
\res_{r,s}^L\Delta_{r,s}(\mu^L_{\lambda},\mu^R_{\lambda})).
\end{array}$$
By the remarks above and Theorem \ref{resrule} this final Hom-space is
isomorphic to
$$\Hom(\Delta_{r-1,s}(\lambda^L-\square,\lambda^R), 
\Delta_{r-1,s}(\mu^L_{\lambda},\mu^R_{\lambda}+\square'))$$
which is non-zero by induction.
\end{proof}

\begin{cor}\label{blocks}
Let $k$ be $\Sigma$-semisimple. 
Two weights $(\lambda^L,\lambda^R)$ and $(\mu^L,\mu^R)$ are in the
same block of $B_{r,s}$ if and only if they are balanced. Each block
contains a unique minimal weight.
\end{cor}
\begin{proof}
One implication was proved in Corollary \ref{oneway}. For the reverse
implication we proceed by induction. If $(\lambda^L,\lambda^R)$
contains a smaller balanced weight $(\mu^L,\mu^R)$ then by Theorem
\ref{maxhom} there exists some
$(\mu^L_{\lambda},\mu^R_{\lambda})\subset (\lambda^L,\lambda^R)$ with
a non-zero homomorphism from $\Delta_{r,s}(\lambda^L,\lambda^R)$ to
$\Delta_{r,s}(\mu^L_{\lambda},\mu^R_{\lambda})$, and hence and
$(\mu^L_{\lambda},\mu^R_{\lambda})$ will lie in the same block. By
induction we also have that $(\lambda^L,\lambda^R)$ and
$(\mu^L,\mu^R)$ lie in the same block. Thus it is enough to show that
there is a unique minimal weight in the set of weights balanced with
$(\lambda^L,\lambda^R)$.

But given two such minimal weights $(\mu^L,\mu^R)$ and
$(\nu^L,\nu^R)$, set $\eta^L=\mu^L\cap\nu^L$ and
$\eta^R=\mu^R\cap\nu^R$. Clearly $(\eta^L,\eta^R)$ is a weight, and
forms a balanced pair with both $(\mu^L,\mu^R)$ and $(\nu^L,\nu^R)$
(and hence with $(\lambda^L,\lambda^R)$). But this contradicts the
minimality of $(\mu^L,\mu^R)$ and $(\nu^L,\nu^R)$.
\end{proof}

In particular, we have now determined the blocks of the walled Brauer
algebra in characteristic zero.

\section{An alcove geometry for the walled Brauer algebra}
\label{geom1}

We would like to have an alcove geometry, coming from some suitable
reflection group, which controls the representation theory of the
walled Brauer algebra in the non-semisimple cases. This is typically
regarded as a Lie theoretic phenomenon, but has been shown to exist
for the Brauer algebra in \cite{cdm2}.

By Theorem \ref{whenss} we may assume that $\delta\in\ZZ$. In this
case we will show that there is such a geometry for the walled Brauer
algebra, associated to $\Sigma_{r+s}$, the Weyl group of type
$A_{r+s}$. However, we will need a new notion of dominant weights (and a
modified group action) to realise this.

Let
$\{\epsilon_{-r},\epsilon_{-(r-1)},\ldots,\epsilon_{-1},\epsilon_1,
\epsilon_2,\ldots,\epsilon_s\}$ be a set of formal symbols. We set
$$X=X_{r,s}=\bigoplus_{i=-r}^{-1}\ZZ\epsilon_i\oplus
\bigoplus_{i=1}^{s}\ZZ\epsilon_i$$
which will be our weight lattice. We will denote an element 
$$\lambda=\lambda_{-r}\epsilon_{-r}+\cdots+\lambda_{-1}\epsilon_{-1}+
\lambda_1\epsilon_1+\cdots+\lambda_s\epsilon_s$$
in $X$ by $(\lambda_{-r},\lambda_{-(r-1)},\ldots,\lambda_{-1};
\lambda_1,\ldots,\lambda_s)$. The set of dominant weights in $X$ is
defined to be
$$X^+=\{\lambda\in X: 0\geq \lambda_{-r}\geq \lambda_{-(r-1)}
\geq\cdots\geq \lambda_{-1}\ \wand \ 
\lambda_1\geq\lambda_2\geq\cdots\geq \lambda_s\geq 0\}.$$
(This is not a standard choice of dominant weights, but will be
justified by our labelling conventions for the walled Brauer algebra.)
Define an inner product on $E=X\otimes_{\ZZ}\RR$ by setting 
$$(\epsilon_i,\epsilon_j)=\delta_{ij}$$
for all nonzero $i,j$ with $-r\leq i,j\leq s$, and extending by
linearity.

There is a root system of type $A$ given by
$$\Phi=\{\pm(\epsilon_i-\epsilon_j):-r\leq i<j\leq s,\quad i,j\neq 0\}.$$
For each root $\beta\in\Phi$ we define a reflection $s_{\beta}$ on $E$
by
$$s_{\beta}(\lambda)=\lambda-\frac{2(\lambda,\beta)}{(\beta,\beta)}
\beta
=(\lambda,\beta)\beta$$ for all
$\lambda\in E$, and let $W$ be the group generated by these
reflections. Then $W$ is just the Weyl group of type $A$, which can be
identified with $\Sigma_{r+s}$.

In algebraic Lie theory it is convenient to shift the action of the
Weyl group on weights relative to some fixed vector $\rho$. While the
same will be true here also, our choice of $\rho$ is rather
different. Fix $\delta\in\ZZ$ and define $\rho=\rho(\delta)\in E$ by
$$\rho=(r,r-1,\ldots,1;\delta,\delta-1,\ldots,\delta-s+1).$$
We consider the dot action of $W$ on $E$ given by 
$$w.\lambda=w(\lambda+\rho)-\rho$$
for all $w\in W$ and $\lambda\in E$. Note that this preserves the
lattice $X$ inside $E$.

A pair of partitions $(\lambda^L,\lambda^R)$ with at most $r$ and $s$ parts
respectively will be identified with a dominant weight $\lambda\in
X^+$ via the map
\begin{equation}\label{ident}
(\lambda^L,\lambda^R)\longmapsto(\bar{\lambda}^L;\lambda^R)=
(-\lambda^L_r,-\lambda^L_{r-1},\ldots -\lambda_1^L;\lambda^R_1,
\lambda^R_2,\ldots,\lambda^R_s).
\end{equation}
 It will also be convenient to have
a graphical representation of elements of $X$. We will represent any
$\lambda\in X$ by a sequence of $r+s$ rows of boxes (to be defined
shortly), with $r$ rows above and $s$ below some fixed horizontal
bar. The $i$th row below this bar will be called row $i$, and the
$i$th row above this bar will be called row $-i$. Columns will be
labelled in increasing order from left to right by elements of $\ZZ$,
and there will be a vertical bar between columns $0$ and $1$. (Note
that there is a column $0$, but no row $0$.) With these conventions,
row $i$ in the representation of $\lambda$ will contain all boxes to
the left of column $\lambda_i$ inclusive.

\begin{figure}[ht]
\includegraphics{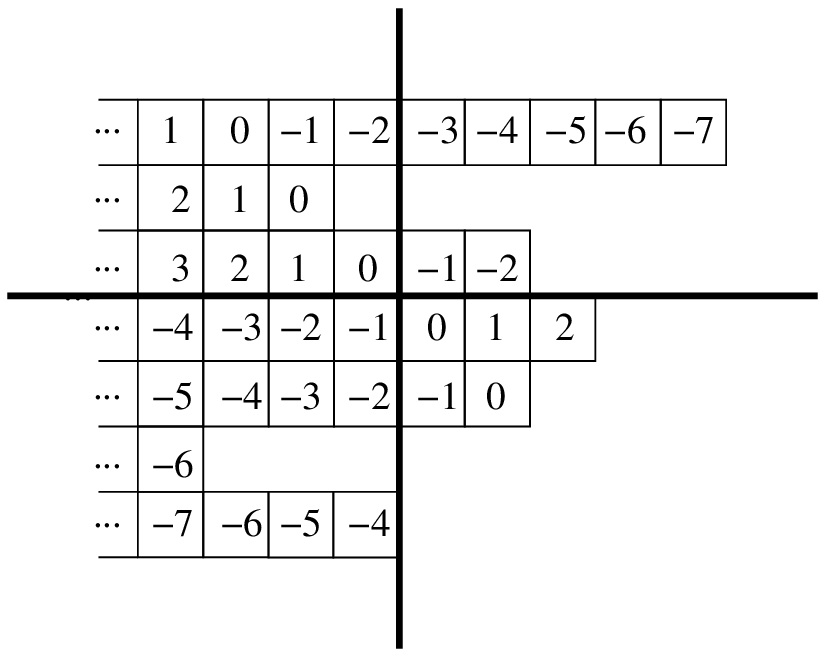}
\caption{The element $(5,-1,2;3,2,-3,0)$ and its associated contents}
\label{egwt}
\end{figure}

We have already defined the content of a box in a pair of partitions
$(\lambda^L,\lambda^R)$. Via the identification in (\ref{ident}) this
corresponds to setting the content of a box in row $i$ and
column $j$ of $\lambda$ to be $j-i$ if $i>0$ and $1+i-j$ if $i<0$.
For example, when $(r,s)=(3,4)$ the element $(5,-1,2;3,2,-3,0)$ (and
the contents of its boxes) is illustrated in Figure \ref{egwt}.

So far our choices of dominant weights, $\rho$, and of contents seems
rather artificial. However, we will see that with these
identifications, the blocks of the walled Brauer algebra for $k$
$\Sigma$-semisimple have a very simple description in terms of
$W$. Before doing this, we will need the following elementary
observation concerning contents of $\delta$-balanced bipartitions. We
denote the number of boxes in a partition $\lambda$ of content $i$ by
$c_i(\lambda)$.

\begin{lem}\label{easycontent}
Two bipartitions $(\mu^L,\mu^R)$ and $(\lambda^L,\lambda^R)$ are
$\delta$-balanced if and only if 
\begin{equation}\label{balcon}
c_i(\lambda^L)-c_i(\mu^L)=c_{-\delta-i}(\lambda^R)-c_{-\delta-i}(\mu^R).
\end{equation}
\end{lem}
\begin{proof}
We begin by noting that the result is obvious if
$(\mu^L,\mu^R)\subseteq (\lambda^L,\lambda^R)$ (or vice versa). In
general, if we have a $\delta$-balanced pair then (\ref{balcon}) is
also clearly satisfied. Thus it is enough to show that the content
condition (\ref{balcon}) implies $\delta$-balanced. 

Suppose that $(\mu^L,\mu^R)$ and $(\lambda^L,\lambda^R)$ satisfy
(\ref{balcon}). Let $\tau^L=\mu^L\cap\lambda^L$ and
$\tau^R=\mu^R\cap\lambda^R$. By the above remarks it is enough to show
that 
\begin{equation}\label{parti}
c_i(\lambda^L)-c_i(\tau^L)=c_{-\delta-i}(\lambda^R)-c_{-\delta-i}(\tau^R)
\end{equation}
and
\begin{equation}\label{partii}
c_i(\mu^L)-c_i(\tau^L)=c_{-\delta-i}(\mu^R)-c_{-\delta-i}(\tau^R).
\end{equation}
Note that $c_i(\tau^L)=\min(c_i(\lambda^L),c_i(\mu^L))$ and
$c_i(\tau^R)=\min(c_i(\lambda^R),c_i(\mu^R))$. Further by (\ref{balcon})
\begin{equation}\label{lammin}
c_i(\tau^L)=c_i(\lambda^L)\quad\mbox{\rm if and only if}\quad
c_{-\delta-i}(\tau^R)=c_{-\delta-i}(\lambda^R).
\end{equation}
Hence exactly one of (\ref{parti}) and (\ref{partii}) becomes the
trivial equality $0=0$.

 Now 
(\ref{balcon}) implies that
$$[c_i(\lambda^L)-c_i(\tau^L)]-[c_i(\mu^L)-c_i(\tau^L)]=
[c_i(\lambda^R)-c_i(\tau^R)]-[c_i(\mu^R)-c_i(\tau^R)].$$
As we already know that one of (\ref{parti}) and (\ref{partii}) holds,
the other is now obvious.
\end{proof}

Lemma \ref{easycontent} will allow us to extend the notion of
$\delta$-balanced to arbitrary pairs of elements in $X$. We cannot
count the number of boxes of a given content in a composition, as
there is no natural point at which to stop including boxes on the left
of the diagram. However we can sensibly extend the notation to
arbitrary differences of compositions as follows. Suppose that
$(\lambda^l;\lambda^r)$ and $(\mu^l;\mu^r)$ are both in $X$ (where we
use lower case superscripts to emphasise that these need not be
bipartitions). We define $c_{i}(\lambda^r-\mu^r)$ (and
$c_{i}(\lambda^l-\mu^l)$) in the following manner. For $1\leq j\leq s$
let $c_{i,j}(\lambda^r-\mu^r)$ be the number of boxes of content $i$
in row $j$ between columns $\min(\lambda^r_j,\mu^r_j)+1$ and
$\max(\lambda^r_j,\mu^r_j)$ inclusive, and let
$\epsilon_{i,j}(\lambda^r-\mu^r)$ be $+1$ if $\lambda^r_j>\mu^r_j$ and
$-1$ otherwise. Then we set
$$c_{i}(\lambda^r-\mu^r)=
\sum_{j=1}^s\epsilon_{i,j}(\lambda^r-\mu^r)c_{i,j}(\lambda^r-\mu^r).$$
Similarly, for $-r\leq j\leq -1$ let $c_{i,j}(\lambda^l-\mu^l)$ be the
number of boxes of content $i$ in row $j$ between columns
$\min(\lambda^l_j,\mu^l_j)+1$ and $\max(\lambda^l_j,\mu^l_j)$
inclusive, and let $\epsilon_{i,j}(\lambda^l-\mu^l)$ be $-1$ if
$\lambda^l_j>\mu^l_j$ and $+1$ otherwise. (Note that this is the
opposite of the previous choice.) Then we set
$$c_{i}(\lambda^l-\mu^l)=
\sum_{j=-r}^{-1}\epsilon_{i,j}(\lambda^l-\mu^l)c_{i,j}(\lambda^l-\mu^l).$$

With these conventions, and using the identification of bipartitions
with $X^+$ in (\ref{ident}), we see that the condition (\ref{balcon})
for bipartitions $(\lambda^L,\lambda^R)$ and $(\mu^L,\mu^R)$ is
equivalent to
\begin{equation}
\label{newbalcon}
c_i(\bar{\lambda}^L-\bar{\mu}^L)=c_{-\delta-i}(\lambda^R-\mu^R).
\end{equation}

We are now in a position to replace our $\delta$-balanced condition
for bipartitions with an orbit condition from the action of the Weyl
group. We proceed in two stages.

\begin{prop}\label{onebway}
Suppose that $(\lambda^L,\lambda^R)$ and $(\mu^L,\mu^R)$ are
bipartions such that 
$$(\bar{\mu}^L;\mu^R)=w.(\bar{\lambda}^L;\lambda^R).$$
Then $(\lambda^L,\lambda^R)$ and $(\mu^L,\mu^R)$ are $\delta$-balanced.
\end{prop}
\begin{proof}
By the above remarks it will be enough to show that (\ref{newbalcon})
holds for any elements of $X$ when $w$ is one of the generators
$s_{\beta}$. Such a generator will only change
$\lambda=(\bar{\lambda}^L;\lambda^R)$ in two rows; there are two cases
to consider depending on whether the two rows are on opposite sides of
the horizontal bar.

From (\ref{newbalcon}) it will be enough to show (i) that in the case
where the two rows are on the same side of the bar the action of
$s_{\beta}$ corresponds to replacing boxes in one row by boxes with
the same content in the other, and (ii) that in the case where the two
rows are on opposite sides of the bar the action of $s_{\beta}$
corresponds to replacing boxes in one row by boxes in the other such
that the two sets can be paired up with contents in each pair summing
to $-\delta$.

Consider a general element $\lambda\in X$. Note that the content of
the last box in row $i>0$ (reading from left to right) equals
$\lambda_i-i$, while the content of the last box in row $-j<0$ is
$-\lambda_{-j}-j+1$. Now consider the action of the generators of $W$
on $\lambda$.

First suppose that $i,j>0$, and consider
$s_{\epsilon_i-\epsilon_j}.\lambda$. We have 
$$\begin{array}{ll}
s_{\epsilon_i-\epsilon_j}.\lambda
&=\lambda
-(\lambda_i-\lambda_j+(\delta-i+1)-(\delta-j+1))(\epsilon_i-\epsilon_j)\\
&=\lambda-(\lambda_i-\lambda_j+j-i)(\epsilon_i-\epsilon_j).
\end{array}$$
If $(\lambda_i-\lambda_j+j-i)>0$ then the action of
$s_{\epsilon_i-\epsilon_j}$ removes $(\lambda_i-\lambda_j+j-i)$ boxes
from row $i$ and adds the same number to row $j$. In row $i$ the boxes
to be removed have contents
$$\lambda_i-i-(\lambda_i-\lambda_j+j-i)+1,\ldots
\lambda_i-i-1,\lambda_i-i$$
and on row $j$ the boxes added have contents
$$\lambda_j-j+1,\lambda_j-j+2,\ldots,\lambda_j-j+(\lambda_i-\lambda_j+j-i).$$
Simplifying both theses expressions for the contents we arrive at the
same list:
$$\lambda_j-j+1,\ldots,\lambda_i-i-i,\lambda_i-i$$
and so the number of boxes of each content has remained unchanged. The
case where $(\lambda_i-\lambda_j+j-i)<0$ is similar.

Next suppose that $i,j>0$, and consider
$s_{\epsilon_{-i}-\epsilon_{-j}}.\lambda$. We have that 
$$s_{\epsilon_{-i}-\epsilon_{-j}}.\lambda
=\lambda-(\lambda_{-i}-\lambda_{-j}+i-j)(\epsilon_{-i}-\epsilon_{-j}).$$
Arguing as above we see that if $(\lambda_{-i}-\lambda_{-j}+i-j)>0$
then this number of boxes are removed from row $-i$ and added to row
$-j$, with both the removed and added sets having contents
$$-\lambda_{-j}-j,\ldots,-\lambda_{-i}-i+1.$$
Thus we again see that the number of boxes of each content remains
unchanged. As above, the 
case where $(\lambda_{-i}-\lambda_{-j}+j-i)<0$ is similar.

Finally, suppose that $i,j>0$, and consider
$s_{\epsilon_{i}-\epsilon_{-j}}.\lambda$. We have that 
$$s_{\epsilon_{i}-\epsilon_{-j}}.\lambda
=\lambda-(\lambda_{i}-\lambda_{-j}
+(\delta-i+1)-j)(\epsilon_{i}-\epsilon_{-j}).$$ Suppose that
$\lambda_{i}-\lambda_{-j} +(\delta-i+1)-j>0$ (the other case is
similar). Then $s_{\epsilon_{i}-\epsilon_{-j}}$ removes
$\lambda_{i}-\lambda_{-j} +(\delta-i+1)-j$ boxes from row $i$ and adds
the same number to row $-j$. In row $-j$ the added contents are
$$-\lambda_{-j}-j,-\lambda_{-j}-j-1,\ldots,
-\lambda_{-j}-j-(\lambda_i-\lambda_{-j}+(\delta-i+1)-j)+1$$
and in row $i$ the removed contents are
$$\lambda_i-i-(\lambda_i-\lambda_{-j}+(\delta-i+1)-j+1,\ldots,
\lambda_i-i-1,\lambda_i-i.$$
Simplifying we see that the added contents are
$$-\lambda_{-j}-j,-\lambda_{-j}-j-1,\ldots,-\lambda_i-\delta+i$$
and the removed contents are
$$\lambda_{-j}-\delta+j,\ldots,\lambda_i-i-1,\lambda_i-i.$$
Comparing the corresponding entries in each of these last two
expressions, we see that in each case they sum to $-\delta$, as required.
\end{proof}

We next consider the reverse implication.

\begin{prop}\label{otherway}
If $(\lambda^L,\lambda^R)$ and $(\mu^L,\mu^R)$
are $\delta$-balanced then there exists $w\in W$ such that 
$$(\bar{\mu}^L;\mu^R)=w.(\bar{\lambda}^L;\lambda^R).$$
\end{prop}
\begin{proof}
It is enough to consider the case where $(\lambda^L,\lambda^R)
\supseteq (\mu^L,\mu^R)$. We will proceed by induction on
$|\lambda^L/\mu^L|=|\lambda^L/\mu^L|$, and write $\lambda$ and $\mu$
for the corresponding elements of $X$.

If $|\lambda^L/\mu^L|=1$ and under our identification the unique box
in $\lambda^L/\mu^L$ (respectively in $\lambda^R/\mu^R$) is in row $-j$
(respectively row $i$) then it is easy to verify that
$$s_{\epsilon_i-\epsilon_{-j}}.\lambda=\mu$$
and so we are done.

Next suppose that  $|\lambda^L/\mu^L|>1$. We define the edge of a skew
partition $\tau$ to be those boxes in $\tau$ such that the box
diagonally below and to the right is not in $\tau$. Let $\epsilon$ be the
box of maximal content in the edge of $\lambda^R/\mu^R$. Suppose that 
the corresponding box in $\lambda$ is in row $i$. By the
$\delta$-balanced condition there is a unique box $\epsilon'$ in the edge
of $\lambda^L/\mu^L$ such that $c(\epsilon)+c(\epsilon')=-\delta$. Suppose
that the corresponding box in $\lambda$ is in row $-j$.

Let $\alpha$ be the first box in row $-j$ \emph{not} in $\lambda$.  As
$\alpha$ is in the edge of $\lambda^L/\mu^L$ we can find a matching
box $\alpha'$ on the edge of $\lambda^R/\mu^R$. Say that $\alpha'$ is
in row $l$; we have that $l\geq i$. This configuration is illustrated
schematically in Figure \ref{strip} with the edge of the two skew
partitions shaded grey, and $\mu$ denoted by the curved lines.

\begin{figure}[ht]
\includegraphics{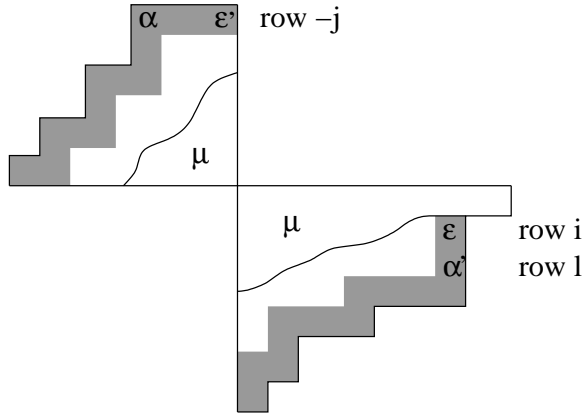}
\caption{Dominant weights $\lambda$ and $\mu$ with associated edge
  strip shaded.}
\label{strip}
\end{figure}

If $i=l$ then define $$\lambda'=s_{\epsilon_i-\epsilon_{-j}}.\lambda.$$
If $i<l$ define
$$\lambda'=(s_{\epsilon_i-\epsilon_{i+1}}\ldots
s_{\epsilon_i-\epsilon_{l-1}}s_{\epsilon_i-\epsilon_{l}}
s_{\epsilon_i-\epsilon_{-j}}).\lambda.$$
In both cases $\lambda'$ is obtained from $\lambda$ by removing the
boxes on the edge of $\lambda^L/\mu^L$ (respectively of
$\lambda^R/\mu^R$) between $\epsilon'$ and $\alpha$ (respectively
between $\epsilon$ and $\alpha'$) inclusive.

If $\lambda'$ is in $X^+$ then we are done by induction. If not, then
either there is a box in the edge of $\lambda^L/\mu^L$ directly below
$\alpha$, or there is a box in the edge of $\lambda^R/\mu^R$ directly
below $\alpha'$. In these cases we repeat the above process replacing
$i$ with $l+1$ and $-j$ with $-j+1$. 

Eventually this process will terminate, as if the whole of the edge of
$\lambda/\mu$ is removed then the result will be in $X^+$. The result
now follows by induction.
\end{proof}

Combining Propositions \ref{onebway} and \ref{otherway} with Corollary
\ref{blocks} we obtain

\begin{cor}\label{geomblocks}
 Let $k$ be $\Sigma$-semisimple with $\delta\in\ZZ$. 
Two weights $\lambda$ and $\mu$ in $X^+$ are in the same block of
  $B_{r,s}$ if and only if $\lambda=w.\mu$ for some $w\in W$.
\end{cor}

\section{Comparing blocks in characteristic zero and characteristic
  $p$}

In this section we will recall some standard results relating blocks
for algebras in different characteristics coming from a common
integral form. This will motivate the results in the following
section.

Let $A_{\ZZ}$ be an integral form giving rise to a
corresponding algebra $A_k$ over a field $k$. Suppose also that we
have a family of integral forms for modules $M_{\ZZ}(\lambda)$ such
that over any given field $k$, all simples can be realised as
quotient modules of the corresponding modules $M_k(\lambda)$ (and are
labelled by the corresponding $\lambda$). We will
have in mind the case when $A_{\ZZ}$ is the integral form for a Brauer
or walled Brauer algebra (where $\delta$ has been specialised to a
fixed integer) and the $M_k(\lambda)$ are cell modules.

Let $1=\sum_ie_i$ be a primitive central idempotent decomposition in
the algebra $A_{\FF_{p}}$ defined over the finite field with $p$
elements. This corresponds to the block decomposition of
$A_{\FF_{p}}$. The idempotent decomposition lifts to a decomposition
over $\ZZ_{p}$ by the lifting theorem (see for example
\cite[Theorem 1.9.4]{bensonI}). These idempotents pass injectively to
$\QQ_p$, and hence to $\CC_p$ which is isomorphic to $\CC$. 

This final decomposition may no longer be primitive, but can be
refined into a primitive decomposition (and hence the blocks over
$\CC$ will in general be smaller). As the labelling scheme for simple
modules in each algebra has been chosen in a consistent manner (via
the integral forms) this proves that if $\lambda$ and $\mu$ are in the
same block over $\CC$, then they are also in the same block over
$\FF_p$.

Combining this with our characteristic zero block results for the
Brauer algebra \cite[Theorem 4.2]{cdm2} and walled Brauer algebra
(Corollary \ref{geomblocks}) we obtain

\begin{prop}\label{evidence}
Let $\lambda$ and $\mu$ be weights for the Brauer (or walled Brauer)
algebra $A(\delta)$ with $\delta\in\ZZ$. Let $W$ be the Weyl group
corresponding to $A$. If there exists $\delta'\in\ZZ$ with
$\delta'\equiv \delta \mod p$ such that $\lambda=w.\mu$ for some $w\in
W$, where the dot action is with respect to $\rho(\delta')$, then
$\lambda$ and $\mu$ are in the same block for $A(\delta)$ over
$\FF_p$.
\end{prop}

\section{A linkage principle in positive characteristic}
\label{geom2}
So far we have given a complete description of the blocks of $B_{r,s}$
when either $k$ is $\Sigma$-semisimple or $\delta\notin \ZZ$. For the
remaining cases we have a necessary conditions for two weights to be in
the same block coming from Corollary \ref{weakdelta}. We will
strengthen this into a linkage principle, using orbits of the affine
Weyl group of type $A$. We will assume throughout this section that
$\delta\in\ZZ$.

Consider the (type $A$) affine Weyl group $W_p$, the group generated
by the affine reflections
$$s_{\beta,rp}(\lambda)=\lambda-((\lambda,\beta)-rp)\beta$$ where
$\beta\in\Phi$ and $r\in\ZZ$. Just as for $W$, this acts on $X$ via
the dot action 
$$w.\lambda=w(\lambda+\rho)-\rho.$$ 

It is easy to verify (as in \cite[Section 5]{cdm2}) that the dot
action of $W_p$ on $X$ is generated by the various dot actions of $W$
with respect to $\rho(\delta+rp)$ for $r\in\ZZ$. Thus from Proposition
\ref{evidence}
we might expect this
affine Weyl group action to control in large part the block structure
in positive characteristic. 
We will see that it does in fact give a
necessary condition for two weights to be in the same block. For this
it will be convenient to have the following simple combinatorial
description of when two elements of $X$ are in the same $W_p$-orbit.

Identify $\Sigma_{r+s}$ as usual with the group of permutations of
$\{-r,\ldots,-2,-1,1,2,\ldots,s\}$. Given an element
$\sigma\in\Sigma_{r+s}$ we define an element $\swap(\sigma)_i$ for
each $-r\leq i\leq s$ with $i\neq 0$ by
$$\swap(\sigma)_i=\left\{\begin{array}{rl}
-1 & \wif i<0\ \wand \sigma(i)>0\\
+1 & \wif i>0\ \wand \sigma(i)<0\\
0 & \otherwise.\end{array}\right.
$$ We also generalise the notion of degree from partitions to elements
of $X$ by setting $$|\lambda|=\sum_{i=-r}^s\lambda_i.$$

\begin{lem}\label{combwp}
Suppose that $\lambda,\mu\in X$. Then $\mu\in W_p.\lambda$ if and only if
$|\lambda|=|\mu|$ and there exists $\sigma\in\Sigma_{r+s}$ such that
for all $-r\leq i\leq s$ with $i\neq 0$ we have
$$\mu_i-i=\lambda_{\sigma(i)}-\sigma(i)+\swap(\sigma)_i(\delta+1)\mod p.$$ 
\end{lem}
\begin{proof}
We have $\mu\in W_p.\lambda$ if and only if 
$$\mu+\rho=w(\lambda+\rho)+p\nu$$ for some $w\in W$ and
$\nu\in\ZZ\Phi$.  Considering each component of $\mu$ in turn and
substituting the corresponding values for $\rho$ we obtain the
congruences given in the lemma. The additional condition that
$|\lambda|=|\mu|$ follows by summing over the expressions for each
$\mu_i$, and using the fact that $\nu\in\ZZ\phi$ implies
that $|\nu|=0$.
\end{proof}

Note that under our correspondence with bipartitions, the condition
that $|\lambda|=|\mu|$ is equivalent to the condition that
$$\lambda^R-\lambda^L=\mu^R-\mu^L$$ which we already know is a
requirement for two bipartitions to be labels of $B_{r,s}$. Thus when
$\lambda$ and $\mu$ both come from bipartitions for $B_{r,s}$ it is
enough to check the congruences in Lemma \ref{combwp} to determine if
they are in the same $W_p$-orbit.

\begin{thm}
Suppose that $\delta\in\ZZ$, and that $(\lambda^L,\lambda^R)$ and
$(\mu^L,\mu^R)$ are bipartitions. If
there exists $M\leq \Delta_{r,s}(\mu^L,\mu^R)$ with 
\begin{equation}\label{ahom}
\Hom(\Delta_{r,s}(\lambda^L,\lambda^R),\Delta_{r,s}(\mu^L,\mu^R)/M)\neq
0
\end{equation}
then $(\bar{\mu}^L;\mu^R)\in W_p.(\bar{\lambda}^L;\lambda^R)$.
\end{thm}
\begin{proof}
By the cellular structure of $B_{r,s}$, if (\ref{ahom}) holds then we
must have $(\lambda^L,\lambda^R)\in\Lambda^a_{r,s}$ and
$(\mu^L,\mu^R)\in\Lambda^b_{r,s}$ for some $b\leq a$. We will proceed by
induction on $r+s$. 

First suppose that $s=0$. Then $B_{r,s}\cong\Sigma_r$, and we have
$$\Delta_{r,s}(\lambda^L,\lambda^R)\cong S^{\lambda^L}\quad\wand\quad
\Delta_{r,s}(\mu^L,\mu^R)\cong S^{\mu^L}.$$ By the block result for
symmetric groups (see for example the formulation in \cite{schurIV})
there exists $v\in W_p^{A_r}$, the affine Weyl group of type $A_r$
such that $\mu^L=v.\lambda^L$. But this affine Weyl group is a
subgroup of $W_p$, and taking $w$ to be the corresponding element in
$W_p$ we have that $\lambda=w.\mu$. A similar argument holds when
$r=0$. 

Now suppose that $\lambda^L$ and $\lambda^R$ are both nonempty
partitions. By localising we may assume that
$(\lambda^L,\lambda^R)\in\Lambda^{r,s}$. Possibly by enlarging $M$, we
may also assume that our non-zero homomorphism kills everything in
$\Delta_{r,s}(\lambda^L,\lambda^R)$ except for one simple in the
head. Such a simple is labelled by some weight
$(\tau^L,\tau^R)\in\Lambda^{r,s}_{reg}$, and by the block result for
symmetric groups, this weight is in the same $W_p$ orbit as
$(\lambda^L,\lambda^R)$. Clearly there will also be a homomorphism
from $\Delta_{r,s}(\tau^L,\tau^R)$ into $\Delta_{r,s}(\mu^L,\mu^R)/M$,
and so we may assume that $
(\lambda^L,\lambda^R)\in\Lambda^{r,s}_{reg}$.

Choose the highest removable box $\Box$ in $\lambda^R$. 
The partition $\lambda^R-\Box$ cannot have a higher addable box of the
same content as $\Box$ by our assumption on
$(\lambda^L,\lambda^R)$. Therefore by the block result for symmetric
groups together with Corollary \ref{indrule} (and Remark \ref{better})
we have a
surjection
$$\ind^R\Delta_{r,s-1}(\lambda^L,\lambda^R-\Box)\too
\Delta_{r,s}(\lambda^L,\lambda^R)\too 0
$$
and so by (\ref{ahom}) we have
$$\Hom(\ind^R\Delta_{r,s-1}(\lambda^L,\lambda^R-\Box),
\Delta_{r,s}(\mu^L,\mu^R)/M)\neq
0.$$ Applying Frobenius reciprocity we see that
$$\Hom(\Delta_{r,s-1}(\lambda^L,\lambda^R-\Box),
\res^R(\Delta_{r,s}(\mu^L,\mu^R)/M))\neq 0.$$
We will set $\lambda'^L=\lambda^L$ and $\lambda'^R=\lambda^R-\Box$.
By Theorem \ref{resrule} we must have either 
\begin{equation}\label{caseone}
\Hom(\Delta_{r,s-1}(\lambda'^L,\lambda'^R),
\Delta_{r,s-1}(\mu^L+\Box',\mu^R)/N)\neq 0
\end{equation}
 for some addable box
$\Box'$ for $\mu^L$ and some $N<\Delta_{r,s-1}(\mu^L+\Box',\mu^R)$ or
\begin{equation}\label{casetwo}
\Hom(\Delta_{r,s-1}(\lambda'^L,\lambda'^R),
\Delta_{r,s}(\mu^L,\mu^R-\Box')/N)\neq 0
\end{equation}
for some removable box 
$\Box'$ in $\mu^R$ and some $N<\Delta_{r,s-1}(\mu^L,\mu^R-\Box')$.

First suppose that we are in the situation in (\ref{caseone}), and set
$(\tau^L,\tau^R)=(\mu^L+\Box',\mu^R)$. The condition in Corollary
\ref{weakdelta} applies both to the weights $(\lambda^L,\lambda^R)$
and $(\mu^L,\mu^R)$ and to the weights $(\lambda'^L,\lambda'^R)$ and
$(\tau^L,\tau^R)$. Comparing the resulting expressions we see that we
must have
$$c(\Box)+c(\Box')+\delta=0\mod p.$$
Also by induction there exists $w\in W_p$ such that $\tau=w.\lambda'$

Alternatively, suppose that we are in the situation in
(\ref{casetwo}), and set $(\tau^L,\tau^R)=(\mu^L,\mu^R-\Box')$. As in
the preceding case, we deduce from Corollary \ref{weakdelta} that
$$c(\Box)=c(\Box')\mod p.$$ Also by induction there exists $w\in W_p$
such that $\tau=w.\lambda'$.

In both cases, suppose that $\Box$ is in row $i$ of $\lambda$ and
$\Box'$ is in row $j$ of $\tau$, and that $\sigma$ is the element of
$W$ corresponding to $w$ as in Lemma \ref{combwp}. Let $\sigma(j)=t$
and $\sigma(u)=j$ for some $t$ and $u$. Define $\sigma'$ by setting
$\sigma'(j)=i$, $\sigma'(u)=t$, and $\sigma'(v)=\sigma(v)$ for all
$v\neq j,u$. It is easy to check that $\sigma'$ satisfies the
conditions in Lemma \ref{combwp}, and hence that $\mu=w'.\lambda$ for
some $w'\in W_p$ as required.
\end{proof}

An immediate consequence of this (and the cellularity of $B_{r,s}$) is

\begin{cor}
Suppose that $\delta\in\ZZ$. Two simples
$L_{r,s}(\lambda^L,\lambda^R)$ and $L_{r,s}(\mu^L,\mu^R)$ are in the
same block only if $(\bar{\mu}^L;\mu^R)\in
W_p.(\bar{\lambda}^L;\lambda^R)$.
\end{cor}

\section{Concluding remarks}

We have given a complete characterisation of the blocks of the walled
Brauer algebra in characteristic zero, and a linkage principle in
characteristic $p$. In general the positive characteristic result
cannot be strengthened to give the full blocks as orbits of the affine
Weyl group, as this is not true for the special case of the ordinary
Brauer algebra (see \cite[Theorem 7.2]{cdm2}).

The geometric description of these results depends on our choice of
embedding of  dominant weights inside a larger weight
space. This has two aspects: the use of \lq negative\rq\ partitions
and the relative positions of the left and right hand components of a
bipartition.

Negative partitions are used so that the natural action of the
symmetric group (which would normally preserve the total number of
boxes) now correspond to adding or removing boxes when acting on rows
from both parts of the bipartition.

The choice of relative positions of the two parts is slightly more
arbitrary. In particular, \cite[Section 4.5]{ddratschur} adopts an
alternative convention of placing the negative partition below the
usual partition. If $\delta> r+s$ then this can be done in such a way
that the Weyl group action we describe corresponds to the standard
choice of $\rho$ from Lie theory (i.e. without a shift by
$\delta$). However this is precisely the case where the block result
is trivial, as all blocks consist of singletons. In general there is
no way to position negative partitions so that $\rho$ is the standard
shift from Lie theory, as for small values of $\delta$ the rows would
have to overlap.

The convention used in this paper also has the advantage that the
embedding of $B_{r,s}$ into $B_{r+t,s+t}$ by globalising is compatible
with the natural embedding of weight spaces for each algebra.

\bibliographystyle{amsalpha} \bibliography{/home/anton/Work/Lib/books}

\newcommand{\etalchar}[1]{$^{#1}$}
\providecommand{\bysame}{\leavevmode\hbox to3em{\hrulefill}\thinspace}
\providecommand{\MR}{\relax\ifhmode\unskip\space\fi MR }
\providecommand{\MRhref}[2]{%
  \href{http://www.ams.org/mathscinet-getitem?mr=#1}{#2}
}
\providecommand{\href}[2]{#2}
\begin{thebibliography}{CMPX06}

\bibitem[BCH{\etalchar{+}}94]{bchlls}
G.~Benkart, M.~Chakrabarti, T.~Halverson, R.~Leduc, C.~Lee, and J.~Stroomer,
  \emph{Tensor product representations of general linear groups and their
  connections with {Brauer} algebras}, J. Algebra \textbf{166} (1994),
  529--567.

\bibitem[Ben91]{bensonI}
D.~J. Benson, \emph{Representations and cohomology {I}}, Cambridge studies in
  advanced mathematics, vol.~30, CUP, 1991.

\bibitem[Bra37]{brauer}
R.~Brauer, \emph{On algebras which are connected with the semisimple continuous
  groups}, Ann. of Math. \textbf{38} (1937), 857--872.

\bibitem[CDM05]{cdm}
A.~G. Cox, M.~{De Visscher}, and P.~P. Martin, \emph{The blocks of the {B}rauer
  algebra in characteristic zero}, preprint, 2005.

\bibitem[CDM06]{cdm2}
\bysame, \emph{A geometric characterisation of the blocks of the {B}rauer
  algebra}, preprint, 2006.

\bibitem[CMPX06]{cmpx}
A.~G. Cox, P.~P. Martin, A.~E. Parker, and C.~Xi, \emph{Representation theory
  of towers of recollement: theory, notes, and examples}, J. Algebra
  \textbf{302} (2006), 340--360.

\bibitem[CR81]{cr1}
C.~W. Curtis and I.~Reiner, \emph{Methods of representation theory}, vol.~1,
  Wiley, 1981.

\bibitem[DD05]{ddratschur}
R.~Dipper and S.~Doty, \emph{The rational {S}chur algebra}, preprint, 2005.

\bibitem[DDH]{ddh}
R.~Dipper, S.~Doty, and J.~Hu, \emph{Brauer's centralizer algebras, symplectic
  {S}chur algebras and {Schur-Weyl} duality}, Trans. AMS, to appear.

\bibitem[Dia88]{diaconis}
P.~Diaconis, \emph{Group representations in probability and statistics},
  Institute of Mathematical Statistics Lecture Notes, Monograph Series 11,
  1988.

\bibitem[Don94]{schurIV}
S.~Donkin, \emph{On {S}chur algebras and related algebras {IV}: {T}he blocks of
  the {S}chur algebras}, J. Algebra \textbf{168} (1994), 400--429.

\bibitem[DWH99]{dhw}
W.~F. Doran, D.~B. Wales, and P.~J. Hanlon, \emph{On the semisimplicity of the
  {B}rauer centralizer algebras}, J. Algebra \textbf{211} (1999), 647--685.

\bibitem[Eny02]{enyang}
J.~Enyang, \emph{Bases of certain algebras associated with quantum groups},
  Ph.D. thesis, Univ. Illinois at Chicago, 2002.

\bibitem[Ful97]{fultab}
W.~Fulton, \emph{Young tableaux}, LMS Student Texts, vol.~35, Cambridge, 1997.

\bibitem[GL96]{gl}
J.~J. Graham and G.~I. Lehrer, \emph{Cellular algebras}, Invent. Math.
  \textbf{123} (1996), 1--34.

\bibitem[GM07]{margreentab}
R.~M. Green and P.~P. Martin, \emph{Constructing cell data for diagram
  algebras}, J. Pure and Applied Algebra \textbf{209} (2007), 551--569.

\bibitem[Gre80]{green}
J.~A. Green, \emph{Polynomial representations of {GL$_n$}}, Lecture Notes in
  Mathematics 830, Springer, 1980.

\bibitem[Hal96]{halvwall}
T.~Halverson, \emph{Characters of the centralizer algebras of mixed tensor
  representations of {$G(r,\mathbb C)$} and the quantum group {$U_q( {\mathfrak
  {gl}}(r,\mathbb C))$}}, Pacific J. of Mathematics \textbf{174} (1996),
  359--410.

\bibitem[HW90]{hw3}
P.~J. Hanlon and D.~B. Wales, \emph{Computing the discriminants of {B}rauer's
  centralizer algebras}, Math. Comp. \textbf{54} (1990), 771--796.

\bibitem[Jam78]{james}
G.~D. James, \emph{The representation theory of the symmetric groups}, Lecture
  Notes in Mathematics 682, Springer, 1978.

\bibitem[JK81]{jk}
G.~D. James and A.~Kerber, \emph{The representation theory of the {S}ymmetric
  group}, Encyclopedia of Mathematics and its Applications, vol.~16,
  Addison-Wesley, 1981.

\bibitem[Koi89]{koikewall}
K.~Koike, \emph{On the decomposition of tensor products of the representations
  of classical groups: by means of universal characters}, Adv. in Math.
  \textbf{74} (1989), 57--86.

\bibitem[KX98]{kxcell}
S.~K\"onig and C.~Xi, \emph{On the structure of cellular algebras}, Algebras
  and modules II (Geiranger, 1996), CMS Conf. Proc., vol.~24, AMS, 1998,
  pp.~365--386.

\bibitem[KX99]{kxinfl}
\bysame, \emph{Cellular algebras: inflations and {M}orita equivalences}, J.
  London Math Soc. \textbf{60} (1999), 700--722.

\bibitem[KX01]{kxbrauer}
\bysame, \emph{A characteristic free approach to {B}rauer algebras}, Trans. AMS
  \textbf{353} (2001), 1489--1505.

\bibitem[MS94]{msblob}
P.~P. Martin and H.~Saleur, \emph{The blob algebra and the periodic
  {T}emperley-{L}ieb algebra}, Lett. Math. Phys. \textbf{30} (1994), 189--206.

\bibitem[Rui05]{ruibrauer}
H.~Rui, \emph{A criterion on the semisimple {B}rauer algebras}, J. Comb. Theory
  Ser. A \textbf{111} (2005), 78--88.

\bibitem[Tur89]{turwall}
V.~Turaev, \emph{Operator invariants of tangles and {$R$-matrices}}, Izvestija
  AN SSSR ser. math. \textbf{53} (1989), 1073--1107, (in Russian).

\bibitem[Wey46]{weyl}
H.~Weyl, \emph{The classical groups, their invariants and representations},
  Princeton University Press, 1946.

\end{thebibliography}
\end{document}